\documentclass[12pt]{article}

\usepackage{amssymb}
\usepackage{amsmath}
\usepackage{graphicx}
\usepackage{color}
 \newtheorem{thm}{Theorem}[section]
 \newtheorem{cor}[thm]{Corollary}
 \newtheorem{lem}[thm]{Lemma}
 \newtheorem{prop}[thm]{Proposition}
 \newtheorem{defn}[thm]{Definition}
 \newtheorem{rem}[thm]{Remark}
 
 \numberwithin{equation}{section}

\newcommand{\tp}{\overset{\!\!\!\!\!\circ}{\sT_+}}
\newcommand{\ind}{\mathrm{ind}\,}
\newcommand{\Span}{\mathrm{span}\,}
\newcommand{\diag}{\mathrm{diag}\,}

 \newcommand{\vp}{\varphi}
 \newcommand{\ve}{\varepsilon}
\newcommand{\D}{\displaystyle}

\newcommand{\smb}{\mbox{\rm smb}\,}
\newcommand{\coker}{\mbox{\rm coker}\,}

\renewcommand{\Im}{\mathrm{Im}\,}

\newcommand{\im}{\mathrm{im}\,}

\newcommand{\nn}{\nonumber}

\newcommand{\re}{\textrm{Re}\,}

 \newcommand{\esssup}{\hbox{\rm ess}\sup_{\!\!\!\!\!\!\!\!\!\!\! t\in {\mathbb{T}}}}
 \newcommand{\rbx}{\hfill{\rule{1.8mm}{1.8mm}}}

\newcommand{\cA}{\mathcal{A}}
\newcommand{\cB}{\mathcal{B}}

\newcommand{\cL}{\mathcal{L}}

\newcommand{\sC}{{\mathbb C}}

\newcommand{\sN}{{\mathbb N}}

\newcommand{\sR}{{\mathbb R}}

\newcommand{\sT}{{\mathbb T}}

\newcommand{\sZ}{{\mathbb Z}}

\begin{document}

%
%
%
%

\vspace*{10mm}

\begin{center}
{\huge Structure of Kernels and Cokernels\\[1ex]
 of Toeplitz plus Hankel
Operators}
\end{center}

\vspace{5mm}

\begin{center}

\textbf{Victor D. Didenko and Bernd Silbermann}\footnote{This
research was supported by the Universiti Brunei Darussalam under
Grant UBD/GSR/S\&T/19.}


\vspace{2mm}


 Universiti Brunei Darussalam,
Bandar Seri Begawan, BE1410  Brunei; diviol@gmail.com

Technische Universit{\"a}t Chemnitz, Fakult{\"a}t f\"ur Mathematik,
09107 Chemnitz, Germany; silbermn@mathematik.tu-chemnitz.de

 \end{center}


  \vspace{10mm}




\begin{abstract}
Toeplitz plus Hankel operators $T(a)+H(b)$, $a,b\in L^\infty$ acting
on the classical Hardy spaces $H^p$, $1<p<\infty$, are studied. If
the generating functions $a$ and $b$ satisfy the so-called matching
condition
 $$
a(t) a(1/t)=b(t) b(1/t),
 $$
an effective description of the structure of the kernel and cokernel
of the corresponding operator is given.  The results depend on the
behaviour of two auxiliary scalar Toeplitz operators, and if the
generating functions $a$ and $b$ are piecewise continuous, more
detailed results are obtained.
\end{abstract}

\section{Introduction\label{s1}}

Let $X$ be a Banach space, and let $\cL(X)$ be the Banach algebra of
all linear continuous operators on $X$. An operator $A\in\cL(X)$ is
called Fredholm if the range $\im A:=\{ Ax:x\in X\}$ of the operator
$A$ is a closed subset of $X$ and the null spaces $\ker A:=\{x\in X:
Ax=0\}$ and $\ker A^*:=\{h\in X^*: A^*h=0\}$ of the operator $A$ and
the adjoint operator $A^*$ are finite-dimensional. For the sake of
convenience, the null space of the adjoint operator $A^*$ is called
the cokernel of $A$ and it is denoted by $\coker A$. Further, if an
operator $A\in\cL(X)$ is Fredholm, then the number
 \begin{equation*}
\kappa:=\dim\ker A-\dim\coker A\, ,
\end{equation*}
where $\dim Y$ denotes the dimension of the linear space $Y$, is
referred to as the index of the operator $A$.

 \sloppy

The present paper is devoted to Toeplitz plus Hankel operators
acting on the classical Hardy spaces. Note that the theory of
Toeplitz and Hankel operators has a long and interesting history and
is distinguished by exiting results and rich connections with many
fields of mathematics, physics, statistical mechanics, and so on
(see, for example, \cite{BS, Pe:2003}). Despite the fact that
Toeplitz operators and Hankel operators are quite different in their
nature, they are closely related to each other. No wonder that
Toeplitz plus Hankel operators $T(a)+H(b)$ have attracted great
attention, as well. Some particular points of interest are Fredholm
properties, index, and even the invertibility. This is caused not
only by theoretical research but also by interesting and challenging
problems which arise in applications and which can be described in
terms of the operators mentioned. For example, Wiener--Hopf plus
Hankel operators occur in scattering theory \cite{MST:1992}, whereas
Toeplitz plus Hankel operators are used in the theory of random
matrix ensembles \cite{BCE:2010, BE:2009}. In the latter case, an
interesting class of Toeplitz plus Hankel operators $T(a)+H(b)$ has
been considered. This class is defined by the condition
 \begin{equation}\label{eq1.1}
a(t) a(1/t)=b(t) b(1/t), \quad t\in\sT\, ,
\end{equation}
where $\sT:=\{t\in\sC: |t|=1\}$ is the unit circle equipped with
counterclockwise orientation.

At present, there is a well-developed Fredholm theory dealing for
Toe\-plitz plus Hankel operators with piecewise continuous
generating functions acting on various Banach and Hilbert spaces.
For more detailed information the reader can consult \cite[Sections
4.95-4.102]{BS}, \cite[Sections 4.5 and 5.7]{RSS2011},
\cite{RS1990}, and \cite{RS:2012}. The latter paper also contains a
transparent index formula for the operators acting on the space
$l^p(\sZ_+)$. The ideas of that work, with a necessary adjustment,
can also be used to obtain similar formulas for operators acting on
Hardy spaces $H^p$. It is worth noting that contrary to scalar
Toeplitz operators, Fredholm Toeplitz plus Hankel operators are not
necessarily one-sided invertible. In a sense, their behaviour is
similar to the behaviour of block Toeplitz operators. Let us recall
that the kernel and cokernel dimension of a block Toeplitz operator
acting on $H^p$ spaces  can be expressed via partial indices of the
Wiener--Hopf factorization of the corresponding generating matrix.
However, as a matter of fact, for an arbitrary matrix-function there
is no efficient procedure in order to obtain its Wiener--Hopf
factorization. On the other hand, for Toeplitz plus Hankel operators
acting on the Hardy spaces $H^p$, there is a general result
\cite{Ehr:2004h, E2004} stating that $T(a)+H(b):H^p\to H^p$, $a,b\in
L^\infty$ is Fredholm if and only if the matrix-function
\begin{equation*}
V(a,b):=\left(%
\begin{array}{cc}
 a- b\widetilde{b}\widetilde{a}^{-1}  & b\widetilde{a}^{-1} \\
 - \widetilde{b}\widetilde{a}^{-1}  & \widetilde{a}^{-1} \\
   \end{array}%
\right),
 \end{equation*}
where $\widetilde{a}=a(1/t), \widetilde{b}=b(1/t)$, admits a certain
type of antisymmetric factorization. Moreover, the defect numbers
$\dim\ker(T(a)+H(b))$ and $\dim\coker(T(a)+H(b))$ of the operator
$T(a)+H(b)$ can be expressed via partial indices of the
antisymmetric factorization of the matrix $V(a,b)$. Nevertheless,
for  arbitrary functions $a,b \in L^\infty$ the chances of finding
such a factorization and the corresponding partial indices are very
slight. However, if the generating functions of the operator
$T(a)+H(b)$ satisfy the additional condition \eqref{eq1.1}, the
matrix $V(a,b)$ takes the form
\begin{equation*}
  V(a,b)=\left(%
\begin{array}{cc}
  0 & d \\
  -c  &  \widetilde{a}^{-1}\\
   \end{array}%
\right) ,
\end{equation*}
where $c:=a b^{-1}= \widetilde{b}\widetilde{a}^{-1}$ and $d:=b
\widetilde{a}^{-1}=a \widetilde{b}^{-1}$. Thus $V(a,b)$ becomes a
triangular matrix, which is better suited for factorization theory.
These ideas have been used while considering Toeplitz plus Hankel
operators of the form $T(a)+H(at^{-1})$ (see \cite{Ehr:2004h}). On
the other hand, the study of the operators $T(a)+H(a)$ in
\cite{BE2004, BE2006} does not involve the factorization of the
matrix $V(a,a)$. Recently, a new method to investigate the operators
$T(a)+H(b):H^p\to H^p$ with piecewise continuous generating
functions $a$ and $b$ satisfying condition \eqref{eq1.1} has been
proposed \cite{BE:2013}. This method is based on the antisymmetric
factorization of the scalar functions $c$ and $d^{-1}$ and it leads
to a complete description of Fredholm properties of the operators
under consideration, including the computation of the corresponding
defect numbers. In particular, it is shown that the operator
$T(a)+H(b)$ is Fredholm if and only if the functions $c$ and
$d^{-1}$ can be represented in a special form. Similar problems have
been studied in \cite{DS:2013d}, but the approach of \cite{DS:2013d}
differs from \cite{BE:2013} and does not employ any factorization
theory in the defect number computation. Still, very often it is not
enough to have an information about Fredholmness and index but more
specific results concerning the kernel and cokernel of the
corresponding operator are required. Thus one of the aims of the
present work is to obtain an effective description of the kernels
and cokernels of the operators in question.

In view of this, let us mention paper \cite{KLR:2009} which is the
culmination of the development aimed at Fredholmness and which
contains a description of kernels and cokernels of singular integral
operators with some Carleman backward shifts.  Although the shifts
in \cite{KLR:2009} are slightly different from that appearing in
Hankel operators, the approach of \cite{KLR:2009} is also based on
Wiener--Hopf factorization of $2\times 2$ matrix functions. It is
also worth noting that the corresponding factorization is assumed to
have factorization factors with entries from $L^\infty$. Such an
assumption implies that the corresponding operator is Fredholm on
any space $L^p$, $1<p<\infty$, if and only if it is Fredholm on one
single space, say on $L^2$. Besides, it turns out that most
piecewise continuous generating functions are not covered by this
method.

The approach used in the present paper is completely different from
both \cite{BE:2013} and \cite{KLR:2009}. For example, we try to
avoid Wiener--Hopf factorization techniques as long as possible.
Thus in Section \ref{s3}, assuming only the right invertibility of
the operator $T(c)$, $c=ab^{-1}$ we give a description of the
kernels of the operators $T(a)\pm H(b)$ for $a,b\in L^\infty$.
Moreover, we show which parts of the kernels of the operators
$T(c)$, $c=ab^{-1}$ and $T(d)$, $d=a \widetilde{b}^{-1}$ make a
contribution to the kernels of the operators $T(a)+H(b)$ and
$T(a)-H(b)$. More precisely, there are decompositions of the kernels
of the operators $T(c)$ and $T(d)$ describing these parts. It is
remarkable that Fredholmness of the operator $T(a)$ plays no role in
our considerations. In Section \ref{s4}, some operators of the form
$T(a)\pm H(t^ka)$, $a\in L^\infty$, $k=-1,0,1$ are considered. It is
shown that for such operators a version of Coburn--Simonenko theorem
holds. In particular, Fredholmness of these operators implies their
one-sided invertibility. Note that the proof of these results does
not use any factorization arguments. Therefore, such an approach can
also be used for similar Toeplitz plus Hankel operators with
discontinuous generating functions acting on $l^p$-spaces, where
factorization technique is not available. However, to tell the
truth, in Sections~\ref{s5} and \ref{s6} we use factorizations of
certain scalar functions. But these factorizations are only used in
order to obtain bases for the subspaces participating in the
decompositions mentioned. Such an approach gives a complete
description for the structure of the null spaces of the operators
$T(a)\pm H(b)$ and allows us to determine their dimensions. For this
purpose, we introduce a characteristic of the factorization for
functions $g\in L^\infty$ satisfying the condition $g
\widetilde{g}=1$ and such that $T(g)$ is Fredholm. It is called
factorization signature and takes values $+1$ and $-1$. Similar
parameters occur in literature from time to time, but it seems that
their importance is not truly appreciated so far. In
Section~\ref{s5} we also provide sufficient conditions which make
possible an easy computation of the factorization signature. For
piecewise continuous generating functions, a very simple formula for
the factorization signature is established in Section \ref{s8}.
Section \ref{s6} is devoted to the description of the kernels and
cokernels of the operators $T(a)+H(b)$ and $T(a)-H(b)$, $a,b\in
L^\infty$ in the case where  operators $T(c)$ and $T(d)$ are
Fredholm. The latter condition implies Fredholmness of the operators
$T(a)\pm H(b)$ and allows one to get their indices. This section
also provides some results for the operators $T(a)\pm H(b)$ under
additional assumption that $\ind T(c)=\pm1$ and the factorization
signature is equal to one. This generalizes the corresponding
results of Section~\ref{s4}. Section~\ref{s7} is specified to the
case of piecewise continuous generating functions.

The present paper has some intersection with the recent paper
\cite{BE:2013} but our approach is entirely different and leads to
more general and more detailed results  with an additional advantage
that the kernel and cokernel of Toeplitz plus Hankel operators are
completely described.

\section{Spaces and operators\label{s2}}

Let us introduce some operators and spaces we need. As usual, let
$L^\infty(\sT)$ stand for the $C^*$-algebra of all essentially
bounded Lebesgue measurable functions on $\sT$, and let
$L^p=L^p(\sT)$, $1\leq p \leq\infty$ denote the Banach space of all
Lebesgue measurable functions $f$ such that
 \begin{align*}
||f||_p:= &\left ( \frac{1}{2\pi}\int_\sT |f(t)|^p \, dt \right
)^{1/p}, \quad
1\leq p <\infty,  \displaybreak[3]\\
 ||f||_\infty:  = & \,\, \esssup |f(t)|,
\end{align*}
is finite. Further, let $H^p=H^p(\sT)$ and
$\overline{H^p}=\overline{H^p(\sT)}$ refer to the Hardy spaces of
all functions $f\in L^p$  the Fourier coefficients
 \begin{equation*}
\widehat{f}_n:=\frac{1}{2\pi}\int_0^{2\pi}f(e^{i\theta})
e^{-in\theta}\, d\theta
 \end{equation*}
 of which vanish for all $n<0$ and $n>0$, respectively. It is a classical
result that for $p\in (1,\infty)$ the Riesz projection $P$ defined
by
 \begin{equation*}
P: \sum_{n=-\infty}^{\infty} \widehat{f}_n e^{in\theta} \mapsto
\sum_{n=0}^{\infty} \widehat{f}_n e^{in\theta},
 \end{equation*}
is bounded on the space $L^p$ and its range is the whole space
$H^p$. The operator $Q:=I-P$,
 \begin{equation*}
Q:\sum_{n=-\infty}^{\infty} \widehat{f}_n e^{in\theta} \mapsto
\sum_{n=-\infty}^{-1} \widehat{f}_n e^{in\theta}
 \end{equation*}
is also a projection and its range is a subspace of the codimension
one in $\overline{H^p}$.

We also consider the flip operator $J:L^p\mapsto L^p$,
 \begin{equation*}
(Jf)(t):=\overline{t}f(\overline{t}), \quad t\in \sT,
 \end{equation*}
where the bar denotes the complex conjugation. Note that the
operator $J$ changes the orientation, satisfies the relations
 \begin{equation*}
J^2=I, \quad  JPJ=Q, \quad JQJ=P,
 \end{equation*}
and for any $a\in L^\infty$,
 \begin{equation*}
JaJ=\widetilde{a}I.
 \end{equation*}
Further, for any $a\in L^\infty$ consider an operator
$T(a):H^p\mapsto H^p$, $1<p<\infty$ defined by
 \begin{equation*}
T(a): f\mapsto Paf.
 \end{equation*}
The operator $T(a)$ is obviously bounded and
 \begin{equation*}
||T(a)||\leq c_p ||a||_\infty,
 \end{equation*}
where $c_p$ is the norm of the Riesz projection on $L^p$. This
operator is called Toeplitz operator generated by the function $a$.
Toeplitz operators with matrix-valued generating functions acting on
$H^p\times H^p$ are defined similarly.

For $a\in L^\infty$,  the Hankel operator $H(a):H^p\mapsto H^p$,
$1<p<\infty$ is defined by
 \begin{equation*}
H(a): f\mapsto PaQJf.
 \end{equation*}
It is clear that this operator is also bounded, that is
\begin{equation}\label{eqB}
||H(a)||\leq c_p ||a||_\infty,
\end{equation}
with the same constant $c_p$ as in \eqref{eqB}. However, in contrast
to Toeplitz operators, the corresponding generating function $a$ is
not uniquely defined by the operator itself. Further, if $a$ belongs
to the space of all continuous functions $C=C(\sT)$, then the Hankel
operator $H(a)$ is compact on the space $H^p$. Moreover, Hankel
operators are never Fredholm, whereas if a Toeplitz operator
$T(a):H^p\to H^p$ is Fredholm, then $T(a)$ is one-sided invertible.

Let $\sZ_+:=\sN\cup \{0\}$.  Recall that in the natural basis
$\{t^n\}_{n\in\sZ_+}$ of the space $H^p$, $1<p<\infty$, Toeplitz and
Hankel operators with generating function $a\in L^\infty$ can be,
respectively, represented as infinite matrices
$(\widehat{a}_{k-j})_{k,j=0}^\infty$ and
$(\widehat{a}_{k+j+1})_{k,j=0}^\infty$, where $\widehat{a}_k$ is the
$k$-th Fourier coefficient of the function $a$. In this form,
Toeplitz and Hankel operators appear on the spaces $l^p(\sZ_+)$ (see
\cite[Section 2.3]{BS}). However, the study of Toeplitz plus Hankel
operators on the space $l^p(\sZ_+)$ is much more difficult since it
is connected with the multiplier problem.

Let us also recall some other results concerning Toeplitz plus
Hankel operators. Thus, if $T(a)+H(b)\in\cL(H^p)$, then the adjoint
operator acts on the space $H^q$, $p^{-1}+q^{-1}=1$ and
\begin{equation}\label{eq2.0}
(T(a)+H(b))^*=T(\overline{a})+H(\widetilde{\overline{b}}).
\end{equation}
Suppose now that $a$ belongs to the group $GL^\infty$ of invertible
elements from $L^\infty$ and
\begin{equation}\label{eq2.1}
   a \widetilde{a}= b \widetilde{b}.
\end{equation}
This relation plays an important role in what follows. It is called
the matching condition, and if $a$ and $b$ satisfy \eqref{eq2.1},
then the duo $(a,b)$ is called the matching pair. For each matching
pair $(a,b)$ one can assign another matching pair $(c,d)$, where
$c=ab^{-1}$ and $d= b \widetilde{a}^{-1}$. Such a pair $(c,d)$ is
called the subordinated pair for $(a,b)$, and it is easily seen that
the functions which constitutes a subordinated pair have a specific
property, namely $c\widetilde{c}=1=d\widetilde{d}$. In passing note
that these functions $c$ and $d$ can also be expressed in the form
\begin{equation*}
c=\widetilde{b}\widetilde{a}^{-1}, \quad d= \widetilde{b}^{-1} a.
\end{equation*}
Besides, if $(c,d)$ is the  subordinated  pair for a matching pair
$(a,b)$, then $(\overline{d},\overline{c})$ is the subordinated pair
for the matching pair $(\overline{a}, \widetilde{\overline{b}})$
defining the adjoint operator (see \eqref{eq2.0}). Further, a
matching pair $(a,b)$ is called Fredholm, if the Toeplitz operators
$T(c)$ and $T(d)$ are Fredholm.

In the following, any function $g\in L^\infty$ satisfying the
condition
\begin{equation*}
g \widetilde{g}=1
\end{equation*}
is called matching function.

 \begin{lem}\label{l3}
If $a,b\in L^\infty$, then the following relations hold:

\begin{enumerate}

\item If $(a,b)$ is a matching pair with the subordinated pair
$(c,d)$, then $(at^{-n},bt^n)$ is a matching pair with the
subordinated pair $(ct^{-2n},d)$.

 \item If $n\in\sN$, then
 \begin{equation}\label{eq2.2}
   T(a)+H(b)=(T(a t^{-n})+H(bt^n)) T(t^n).
\end{equation}
\end{enumerate}
\end{lem}

\textbf{Proof.} Assertion (i) can be verified straightforward, and
equation \eqref{eq2.2} is a consequence of the well--known
identities
  \begin{equation*}
\begin{aligned}
   T(a_1a_2)=T(a_1)T(a_2)+H(a_1)H(\widetilde{a_2}), \\
    H(a_1a_2)=T(a_1)H(a_2)+H(a_1)T(\widetilde{a_2}),
\end{aligned}
\end{equation*}
and the relation $H(t^n)T(t^n)=0$.
 \rbx

Recall that $a$ is assumed to be invertible in $L^\infty$ and let us
point out that this is always the case when the operator $T(a)+H(b)$
is Fredholm \cite{BE:2013, E2004}. On the other hand, the
invertibility of $a\in L^\infty$ does not automatically implies the
Fredholmness either of the operators $T(a)+H(b)$ or $T(a)$.
Nevertheless, Fredholm properties of Toeplitz operators $T(a)$ can
be described by using the Wiener--Hopf factorization of the
generating function $a$. Assume that $p>1$, $q>1$ are real numbers
such that $p^{-1}+q^{-1}=1$.
 \begin{defn}
We say that a function $a\in L^\infty$ admits a week Wiener--Hopf
factorization in $H^p$, if it can be represented in the form
 \begin{equation}\label{fac}
    a=a_- t^n a_+,
\end{equation}
where $n\in \sZ$, $a_+\in H^q$, $a_+^{-1}\in H^p$, $a_-\in
\overline{H^p}$, $a_-^{-1}\in \overline{H^q}$, and $a_-(\infty)=1$.
 \end{defn}
It is well-known that the weak Wiener--Hopf factorization of a
function $a$ is unique, if it exists.  The functions $a_-$ and $a_+$
are called the factorization factors, and the number $n$ is the
factorization index. If $a\in L^\infty$ and the operator $T(a)$ is
Fredholm, then the function $a$ admits the weak Wiener--Hopf
factorization with $n=-\ind T(a)$ \cite{BS, LiSp:1987}. Moreover, in
this case, the factorization factors possess an additional
property--viz. the linear operator $a_+^{-1}Pa_-^{-1}I$ defined on
$\Span \{t^k:k\in\sZ_+\}$ can be boundedly extended on the whole
space $H^p$. Throughout this paper, such a kind of weak Wiener--Hopf
factorization in $H^p$ is called simply Wiener--Hopf factorization
in $H^p$. Note that if $h$ is a polynomial, then the element
$Pa_-^{-1} h$ is also a polynomial. The following result is
well--known.
 \begin{thm}[see {\cite{BS}}]\label{t1}
If $a\in L^\infty$, then Toeplitz operator $T(a):H^p\to H^p$,
$1<p<\infty$ is Fredholm and $\ind T(a)=-n$ if and only if the
generating function $a$ admits the Wiener--Hopf factorization
\eqref{fac} in $H^p$.
 \end{thm}
Recall that Fredholmness of a Toeplitz operator depends on the space
where this operator acts. The reader can consult \cite{BS,
LiSp:1987} for more details. Another important result is given by
Coburn--Simonenko theorem stating that if $a$ is the non-zero
element in $L^\infty$, then the kernel or cokernel of the operator
$T(a)$ is trivial. Block Toeplitz operators does not possess such a
property and this causes serious difficulties. For Toeplitz plus
Hankel operators situation is similar to block Toeplitz operators.
Therefore, the determination of the numbers $\dim\ker (T(a)+H(b))$
and $\dim\coker (T(a)+H(b))$ becomes a challenging problem, and one
of the aims of this work is to find them for some classes of the
operators $T(a)+H(b)$. In conclusion of this section, let us recall
that one-sided inverses of a Fredholm scalar Toeplitz operator
$T(a)$ can be effectively derived. Let $n$ be the factorization
index. If $n\geq 0$, then $T(a)$ is left--invertible and the
operator $T(t^{-n})T^{-1}(a_0)$, where $a_0:=a t^{-n}$, is one of
the left--inverses for $T(a)$. On the other hand, if $n\leq0$, then
$T(a)$ is right--invertible. For the sake of convenience, in this
paper the notation $T_r^{-1}(a)$ always means the operator
$T^{-1}(a_0)T(t^{-n})$, which is one of right inverses for the
operator $T(a)$. Besides, for $n>0$ the kernel of the operator
$T(t^{-n})$ is the linear span of the monomials $1,t, \cdots,
t^{n-1}$, i.e. $ \ker T(t^{-n})=\Span \{1,t,\cdots, t^{n-1}\}$.

If $T(a)$ is right--invertible and $\dim\ker T(a)=\infty$, then
$T_r^{-1}(a)$ denotes one of right inverses of $T(a)$.

 \section{Kernels of Toeplitz plus Hankel operators \label{s3}}

Let $a,b\in L^\infty$. On the space $H^p$, $1<p<\infty$ consider
Toeplitz plus Hankel operators $T(a)+H(b)$ and $T(a)-H(b)$. As was
already mentioned, in the present paper, the kernel spaces of these
operators are studied under the additional condition \eqref{eq2.1},
connecting the generating functions $a$ and $b$.  This condition
presents a unique possibility to obtain an effective description for
the kernels of the operators $T(a)\pm H(b)$. Nevertheless, let us
start with an auxiliary result for Toeplitz plus Hankel operators
with arbitrary generating functions $a,b \in L^\infty$. By
$V=V(a,b)$ we denote the matrix
\begin{equation*}
V(a,b):=
\left(%
\begin{array}{cc}
 a- b\widetilde{b}\widetilde{a}^{-1} & d \\
 -c   &  \widetilde{a}^{-1}\\
   \end{array}%
\right),
\end{equation*}
where  $c:=\widetilde{b}\widetilde{a}^{-1}$,
$d:=b\widetilde{a}^{-1}$. As usual, it is assumed that the element
$a\in L^\infty(\sT)$ is invertible. In particular, this condition is
satisfied if at least one of the operators $T(a)+H(b)$ or
$T(a)-H(b)$ is Fredholm \cite{BE:2013, E2004}.

Our first concern here is to introduce a formula which is known in
principle -- viz.
\begin{equation}\label{eqTV}
 \left(%
\begin{array}{cc}
  T(a)+H(b)+Q &  0\\
   0 &  T(a)-H(b)+Q\\
   \end{array}%
\right) =
 A(T(V(a,b))+\diag (Q,Q))B\,,
\end{equation}
where $A,B: L^p\times L^p\to L^p\times L^p$ are invertible
operators,
\begin{equation*}
B= \left(%
\begin{array}{cc}
  I &  0\\
  \widetilde{b}I  &  \widetilde{a}I\\
   \end{array}%
\right)
\left(%
\begin{array}{cc}
  I & I\\
  J  &  -J\\
   \end{array}%
\right),
\end{equation*}
and $A$ is also known but its concrete form is not important right
now. Note that in the following we frequently denote the operator
$aI$ of multiplication by the function $a\in L^\infty$ simply by
$a$.

In order to establish formula \eqref{eqTV}, let us recall some well
known operator identities and relations connecting Toep\-litz plus
Hankel operators with block Toeplitz operators. Thus if $a$ belongs
to a unital Banach algebra with the identity $e$ and if $p$ is an
idempotent from the same algebra, i.e. $p^2=p$, then
\begin{equation}\label{eqn1}
\begin{aligned}
   pa+q  &= (e+paq)(pap+q), \\
pap+q&=(e-paq)(pa+q) ,
\end{aligned}
\end{equation}
where $q=e-p$, and $e-paq$ is the inverse for the element $e+paq$.

On the other hand, one can write (see, for example, \cite{DS:2012a,
Kr:1987})

\begin{multline}\label{eqn3}
 \frac{1}{2} \left(%
\begin{array}{cc}
  I & I \\
  J& -J \\
\end{array}%
\right)
\left(%
\begin{array}{cc}
  T(a)+H(b)+Q & 0 \\
  0 & T(a)-H(b)+Q \\
\end{array}%
\right)
 \left(%
\begin{array}{cc}
  I & J \\
 I & -J \\
\end{array}
\right)  \\
=\diag(P,Q)\left(%
\begin{array}{cc}
  a & b \\
  \widetilde{b}  & \widetilde{a} \\
   \end{array}%
\right)
  \diag(P,Q) +(\diag(I,I)-\diag(P,Q)).
  \end{multline}

Setting $p:=\diag(P,Q)$, $e=\diag(I,I)$ and using \eqref{eqn1} one
can rewrite the relation \eqref{eqn3} as
 \begin{multline}\label{idt}
 \frac{1}{2} \left(%
\begin{array}{cc}
  I & I \\
  J& -J \\
\end{array}%
\right)
\left(%
\begin{array}{cc}
  T(a)+H(b)+Q& 0 \\
  0 & T(a)-H(b)+Q \\
\end{array}%
\right)
 \left(%
\begin{array}{cc}
  I & J \\
 I & -J \\
\end{array}
\right)  \displaybreak[3]\\
=p\left(%
\begin{array}{cc}
  a & b \\
  \widetilde{b}  & \widetilde{a} \\
   \end{array}%
\right)
  p +q  \displaybreak[3]\\
=\left (e-p\left(%
\begin{array}{cc}
  a & b \\
  \widetilde{b}  & \widetilde{a} \\
   \end{array}%
\right)q \right ) \left(e+p\left(%
\begin{array}{cc}
  a & b \\
  \widetilde{b}  & \widetilde{a} \\
   \end{array}%
\right)q \right)
  \left ( p\left(%
\begin{array}{cc}
  a & b \\
  \widetilde{b}  & \widetilde{a} \\
   \end{array}%
\right)
  p +q \right )\\
=\left (e\!-\!p\left(%
\begin{array}{cc}
  a & b \\
  \widetilde{b}  & \widetilde{a} \\
   \end{array}%
\right)q \right )\left ( p\left(%
\begin{array}{cc}
  a & b \\
  \widetilde{b}  & \widetilde{a} \\
   \end{array}%
\right)
   +q \right )\displaybreak[3]\\
= \!\left (e\!-\!p\left(%
\begin{array}{cc}
  a & b \\
  \widetilde{b}  & \widetilde{a} \\
   \end{array}%
\right)q \right ) \left (%
\diag(P,P)\left(%
\begin{array}{cc}
  a & b \\
 0  & 1 \\
   \end{array}%
\right) +
 \diag(Q,Q)\left(%
\begin{array}{cc}
  1 & 0 \\
 \widetilde{b} & \widetilde{a} \\
   \end{array}%
\right)
 \right )\\
 =\left (e\!-\!p\left(\!\!%
\begin{array}{cc}
  a & b \\
  \widetilde{b}  & \widetilde{a} \\
   \end{array}%
\!\!\right )q \right ) \left (%
\diag(P,P)\left(\!\!\!\!%
\begin{array}{c@{\hspace{1mm}}c}
  a- b\widetilde{b}\widetilde{a}^{-1} & b\widetilde{a}^{-1} \\
 -\widetilde{b}\widetilde{a}^{-1}  & \widetilde{a}^{-1} \\
   \end{array}%
\!\!\!\right)\! +\!
 \diag(Q,Q)
 \right )\left(\!\!%
\begin{array}{cc}
  1 & 0 \\
 \widetilde{b} & \widetilde{a}\\
   \end{array}%
\!\!\right)  \displaybreak[3]\\[1ex]
=\!\left (e\!-\!p\left(%
\begin{array}{cc}
  a & b \\
  \widetilde{b}  & \widetilde{a}  \\
   \end{array}%
\right)q \right ) \left (%
\diag(P,P)\left(\!\!%
\begin{array}{cc}
  a- b\widetilde{b}\widetilde{a}^{-1} & d \\
 -c  & \widetilde{a}^{-1} \\
   \end{array}%
\!\!\right)\!+\!
 \diag(Q,Q)
 \right )\left(\!%
\begin{array}{cc}
  1 & 0 \\
 \widetilde{b} & \widetilde{a} \\
   \end{array}%
\!\right) \displaybreak[3]\\[1ex]
=\!\left (e\!-\!p\left(%
\begin{array}{cc}
  a & b \\
  \widetilde{b}  & \widetilde{a} \\
   \end{array}%
\right)q \right ) \left (
 \diag(I,I) \!+ \! \diag(P,P)
 \left(\!\!%
\begin{array}{cc}
 a- b\widetilde{b}\widetilde{a}^{-1} & d \\
 -c  & \widetilde{a}^{-1} \\
   \end{array}%
\!\!\right)
 \diag(Q,Q)
 \right ) \displaybreak[3] \\[1ex]
 \times
 \left (%
\diag(P,P)\left(%
\begin{array}{c@{\hspace{-0.1mm}}c}
  a- b\widetilde{b}\widetilde{a}^{-1} & d \\
 -c  & \widetilde{a}^{-1} \\
   \end{array}%
\right)\diag(P,P) +
 \diag(Q,Q)
 \right )
 \left(%
\begin{array}{cc}
  1 & 0 \\
 \widetilde{b} & \widetilde{a} \\
   \end{array}%
\right).
 \end{multline}
In the relation \eqref{idt}, all the operators connecting the
operator $\diag(T(a)+H(b)+Q,T(a)-H(b)+Q)$ and the operator
 \begin{equation*}
\widetilde{\cA}:=\diag(P,P)\left(%
\begin{array}{cc}
  a- b\widetilde{b}\widetilde{a}^{-1} & d \\
 -c  & \widetilde{a}^{-1} \\
   \end{array}%
\right)\diag(P,P) +
 \diag(Q,Q)
  \end{equation*}
are invertible and
\begin{equation*}
 \left(%
\begin{array}{cc}
  1 & 0 \\
 \widetilde{b} & \widetilde{a} \\
   \end{array}%
    \right )^{-1}=
\left(%
\begin{array}{cc}
  1 & 0 \\
- \widetilde{b}\widetilde{a}^{-1} & \widetilde{a}^{-1} \\
   \end{array}%
    \right ).
\end{equation*}
Thus the representation \eqref{eqTV} is established.

 \begin{rem}\label{r2}
Relation \eqref{eqTV} indicates that the operators
$\diag(T(a)+H(b),T(a)-H(b))$ and $T(V(a,b))$ are simultaneously
Fredholm. However, this conclusion is not always true for the
operators $T(a)+H(b)$ and $T(a)-H(b)$ themselves. Examples can be
already found among the operators $I\pm H(b)$ with piecewise
continuous generating functions (see, for example, \cite{DS:2012a}).
Even if both operators $T(a)+H(b)$ and $T(a)-H(b)$ are Fredholm,
they can have different indices. Thus, in general, the use of
relation \eqref{eqTV} in the study of Fredholm properties of
Toeplitz plus Hankel operators is limited. Nevertheless, this
relation is still very helpful in the investigation of the kernels
of Toeplitz plus Hankel operators.
 \end{rem}

   \begin{lem}\label{l1}
Assume that $a,b\in L^\infty$, $a\in GL^\infty$, and the operators
$T(a)\pm H(b)$ are considered on the space $H^p$, $1<p<\infty$. Then
  \begin{itemize}
    \item If $(\vp,\psi)^T\in \ker T(V(a,b))$, then
  \begin{align}\label{eqKer1}
     (\Phi, \Psi)^T= &
      (\vp-JQc\vp+JQ\widetilde{a}^{-1}\psi,
      \vp+JQc\vp-JQ\widetilde{a}^{-1}\psi)^T
      \\
      &\in
     \ker \diag(T(a)+H(b),T(a)-H(b)) &\nn
\end{align}

    \item If $(\Phi, \Psi)^T\in \ker\diag (T(a)+H(b),T(a)-H(b)
 )$, then
\begin{equation}\label{eqKer2}
(\Phi+\Psi, P(\widetilde{b}(\Phi + \Psi)
+\widetilde{a}JP(\Phi-\Psi))^T\in \ker T(V(a,b)).
\end{equation}
    \end{itemize}
    Moreover, the operators
 \begin{align*}
    E_1: \ker T(V(a,b)) \to \ker\diag(T(a)+H(b),T(a)-H(b)),\\[1ex]
    E_2: \ker\diag(T(a)+H(b),T(a)-H(b)) \to \ker T(V(a,b)),
\end{align*}
defined, respectively, by the relations \eqref{eqKer1} and
\eqref{eqKer2} are mutually inverses to each other.
 \end{lem}

 \textbf{Proof.}
Consider the representation \eqref{eqTV} and note that both
operators $A$ and  $B$ are invertible on the space  $L^p(\sR)\times
L^p(\sR)$.  Therefore, relation \eqref{eqTV} implies that for any
$(\vp,\psi)^T\in \ker T(V(a,b)))$, the element
$B^{-1}((\vp,\psi)^T)$ belongs to the set
 \begin{equation*}
 \ker\diag(T(a)+H(b)+Q,T(a)-H(b)+Q)\!\!=\! \ker\diag(T(a)+H(b),T(a)-H(b)).
 \end{equation*}
 Hence
  \begin{equation*}
\diag (P,P)B^{-1}((\vp,\psi)^T)=B^{-1}((\vp,\psi)^T).
  \end{equation*}
Computing the left-hand side of the last equation, one obtains the
relation \eqref{eqKer1}. Analogously, if $(\Phi, \Psi)^T\in \ker
\diag(T(a)+H(b),T(a)-H(b))$, then $B((\Phi, \Psi)^T) \in \ker
T(V(a,b))$ and
 \begin{equation*}
\diag (P,P) B((\Phi, \Psi)^T)=B((\Phi, \Psi)^T),
 \end{equation*}
so the representation \eqref{eqKer2} follows.

Now let $(\vp,\psi)$ and $(\Phi,\Psi)$ be as above. Then
 \begin{align*}
   \diag (P,P) B \, \diag (P,P) B^{-1}((\vp,\psi)^T)& = B
   B^{-1}((\vp,\psi)^T),\\
 \intertext{and}
 \diag (P,P) B^{-1}\, \diag (P,P) B((\Phi,\Psi)^T)& = B^{-1}
   B((\Phi,\Psi)^T),
\end{align*}
which completes the proof.

  \rbx

From now on we also assume that the generating functions $a,b\in
L^\infty$ satisfy matching conditions \eqref{eq2.1}. If this is the
case, the kernels of Toeplitz plus Hankel operators $T(a)+H(b)$ and
$T(a)-H(b)$ can be studied in more detail.

Note that if $(a,b)$ is a matching pair, then the matrix--function
$V(a,b)$ takes the form
  \begin{equation*}
V(a,b)=\left(%
\begin{array}{cc}
 0  & d \\
  -c  & \widetilde{a}^{-1} \\
   \end{array}%
\right).
\end{equation*}
where $(c,d)$ is the corresponding subordinated pair. In addition,
we also have a useful representation--viz.
 \begin{align}
T(V(a,b))& =\left(%
\begin{array}{cc}
 0  & T(d) \\
  -T(c)  & T(\widetilde{a}^{-1}) \\
   \end{array}%
\right)   \nn\\[1ex]
 &= \left(%
\begin{array}{cc}
 -T(d)  & 0 \\
 0  & I \\
   \end{array}%
\right)
\left(%
\begin{array}{cc}
 0  & -I \\
  I  & T(\widetilde{a}^{-1}) \\
   \end{array}%
\right)
\left(%
\begin{array}{cc}
 -T(c)  &  0\\
  0  & I \\
   \end{array}%
\right),  \label{eq3}
\end{align}
where the operator
 \begin{equation*}
   D:=\left(%
\begin{array}{cc}
 0  & -I \\
  I  & T(\widetilde{a}^{-1}) \\
   \end{array}%
\right)
\end{equation*}
in the right--hand side of \eqref{eq3} is invertible and
 \begin{equation*}
   D^{-1} =\left(%
\begin{array}{cc}
 T(\widetilde{a}^{-1})  & I \\
 - I  & 0 \\
   \end{array}%
\right).
\end{equation*}

Now we can establish a representation for the kernel of the block
Toeplitz operator $T(V(a,b))$.
 \begin{prop}\label{p3.4}
Let $(a,b)\in L^\infty\times L^\infty$ be a matching pair such that
the operator $T(c)$, $c=ab^{-1}$, is invertible from the right. Then
\begin{equation*}
\ker T(V(a,b))=\Omega(c) \dotplus \widehat{\Omega}(d)
\end{equation*}
 where
 \begin{align*}
\Omega(c) &:=\left \{   (\vp,0)^T:\vp\in \ker T(c)\right\},\\
\widehat{\Omega}(d) &:=\left \{ (T_r^{-1}(c)
T(\widetilde{a}^{-1})s,s)^T:s\in \ker T(d) \right \}.
\end{align*}
 \end{prop}
 \textbf{Proof.}
It is easily seen that $\Omega(c)$ and $\widehat{\Omega}(d)$ are
closed subspaces of the kernel of the operator $T(V(a,b))$. It is
also clear that the intersection of this two subspaces consists of
the zero vector only.

 If $(y_1,y_2)^T\in\ker T(V(a,b))$, then
\begin{equation*}
T(d) y_2=0, \text{ and } T(c)y_1=T(\widetilde{a}^{-1})y_2.
\end{equation*}
 Since $T_r^{-1}(c)$ is left-invertible, the space $H^p$ can be
 represented as the direct sum of the closed subspaces $\ker T(c)$
 and $\im T_r^{-1}(c)$, i.e.
 \begin{equation*}
H^p=\ker T(c)\dotplus \im T_r^{-1}(c).
 \end{equation*}
Correspondingly, the element $y_1$ can be written as
$y_1=y_{10}+y_{11}$, where $y_{10}\in\ker T(c)$ and $y_{11}\in \im
T_r^{-1}(c)$. Moreover, there is a unique vector $y_3\in H^p$ such
that $y_{11}=T_r^{-1}(c)y_3$, so we get
\begin{equation*}
T(c) y_1= T(c)(T_r^{-1}(c)y_3 +y_{10})=y_3=
T(\widetilde{a}^{-1})y_2.
\end{equation*}
It implies that
\begin{equation*}
y_1=T_r^{-1}(c)T(\widetilde{a}^{-1})y_2+y_{10},
\end{equation*}
 and we get
 \begin{equation*}
(y_1,y_2)^T=(T_r^{-1}(c)
T(\widetilde{a}^{-1})y_2,y_2)^T+(y_{10},0)^T,
\end{equation*}
with $(T_r^{-1}(c) T(\widetilde{a}^{-1})y_2,y_2)^T\in
\widehat{\Omega}(d)$ and $(y_{10},0)^T\in \Omega(c)$, which
completes the proof. \rbx

We already know that if an element $\vp\neq0$ belongs to the kernel
of the operator $T(c)$, then $(\vp,0)^T \in \ker T(V(a,b))$ and
Lemma \ref{l1} implies that
 \begin{align*}
 \vp-JQcP\vp \in \ker(T(a)+H(b)),\\
 \vp+JQcP\vp\in \ker(T(a)-H(b)).
\end{align*}
However, it is remarkable that the functions $ \vp-JQcP\vp$ and
$\vp+JQcP\vp$ belong to the kernel of the operator $T(c)$ as well.
 \begin{prop}\label{p1}
 Let $g\in L^\infty$ satisfy the relation $g \widetilde{g}=1$. Then
 \begin{enumerate}
\item If $f\in \ker T(g)$, then $JQgPf\in \ker T(g)$.
\item If $f\in \ker T(g)$, then $(JQgP)^2f=f$.
 \end{enumerate}
 \end{prop}

\textbf{Proof.}
 If $g \widetilde{g}=1$ and $f\in \ker T(g)$, then
  \begin{equation*}
  T(g)(JQgPf)=PgPJQgPf
  =JQ\widetilde{g}QgPf=JQ\widetilde{g}gPf-JQ\widetilde{g}PgPf=0,
\end{equation*}
and assertion (i) follows. On the other hand, for any $f\in \ker
T(g)$ one has
\begin{align*}
(JQgP)^2f & =JQgPJQgP f= P\widetilde{g}QgPf&\\
&=P\widetilde{g}gPf- P\widetilde{g}PgPf =f-P\widetilde{g}T(g)f=f,
\end{align*}
and we are done. \rbx

Consider now the operator $\mathbf{P}_g:=JQgP\left |_{\ker
T(g)}\right .$. By Proposition \ref{p1}, one has $\mathbf{P}_g:\ker
T(g)\to \ker T(g)$ and  $\mathbf{P}_g^2=I$. Therefore,  on the space
$\ker T(g)$ the operators $\mathbf{P}_g^{\pm}:=(1/2)(I\pm
\mathbf{P}_g)$ are complimentary projections, so they generate a
decomposition of $\ker T(g)$.

 \begin{cor}\label{c1}
Let $(c,d)$ be the subordinated pair for a matching pair $(a,b)\in
L^\infty\times L^\infty$. Then the following relations
 \begin{align*}
   &\ker T(c)=\im \mathbf{P}_c^{-}\dotplus\im\mathbf{P}_c^{+},\\
  & \im \mathbf{P}_c^{-}\subset \ker(T(a) +H(b)) ,\\
   & \im \mathbf{P}_c^{+}\subset \ker(T(a)-H(b)),
\end{align*}
hold.
 \end{cor}
Corollary \ref{c1} shows the influence of the operator $T(c)$ on the
kernels of the operators $T(a)+H(b)$ and $T(a)-H(b)$. Let us now
clarify the role of another operator--viz. the operator $T(d)$, in
the kernel structure of the corresponding operators $T(a)\pm H(b)$.
This problem is more involved. Assume additionally that the operator
$T(c)$ is invertible from the right. If $s\in \ker T(d)$, then the
element $(T_r^{-1}(c) T(\widetilde{a}^{-1})s,s)^T \in \ker
T(V(a,b))$. By Lemma \ref{l1}, the element
 \begin{equation}\label{eqphi}
   \vp_{\pm}(s):=T_r^{-1}(c)T(\widetilde{a}^{-1}) s \mp JQcP
   T_r^{-1}(c)T(\widetilde{a}^{-1})s \pm JQ \widetilde{a}^{-1} s
\end{equation}
belongs to the null space  $\ker(T(a)\pm H(b))$ of the corresponding
operator $T(a)\pm H(b)$.

 \begin{lem}\label{l2}
The map $s\mapsto \vp_{\pm}(s)$ is a one-to-one function from the
space $\im \mathbf{P}_d^{\pm}$  to the space $\ker(T(a)\pm H(b))$.
 \end{lem}
  \textbf{Proof.}
Assume that $s\in \ker T(d)$. If we show that the operator
$(1/2)(P\widetilde{b}P+P \widetilde{a}JP)$ sends $\vp_{+}(s)$ into
$\mathbf{P}_d^{+}s$ and the operator $(1/2)(P\widetilde{b}P-P
\widetilde{a}JP)$ sends $\vp_{-}(s)$ into $\mathbf{P}_d^{-}s$, then
Lemma \ref{l2} will follow. Consider, for example, the first case.
Thus one has
\begin{equation*}
(P\widetilde{b}P+P \widetilde{a}JP)\vp_+(s)=I_1+I_2+\cdots+I_6,
\end{equation*}
where
 \begin{align*}
  I_1&=P\widetilde{b}P T_r^{-1}(c)T(\widetilde{a}^{-1}) s ,\\
I_2&=-P\widetilde{b}P JQcP T_r^{-1}(c)T(\widetilde{a}^{-1}) s ,\\
I_3&=P\widetilde{b}P JQ\widetilde{a}^{-1}P s ,\\
I_4&=P\widetilde{a}JP T_r^{-1}(c)T(\widetilde{a}^{-1}) s \\
I_5&=-P\widetilde{a}J P JQcP T_r^{-1}(c)T(\widetilde{a}^{-1}) s ,\\
I_6&=P\widetilde{a} JPJ\widetilde{a}^{-1}P s.
\end{align*}
Taking into account that $\widetilde{b}=c\widetilde{a}$, $bc=a$,
$JPJ=Q$, $J\widetilde{a}J=a$ and $P+Q=I$, we get
 \begin{align*}
   I_1+I_5&=P\widetilde{a}c P T_r^{-1}(c)T(\widetilde{a}^{-1}) s-
   P\widetilde{a}QcP T_r^{-1}(c)T(\widetilde{a}^{-1}) s \\
    &=P\widetilde{a}P c P T_r^{-1}(c)T(\widetilde{a}^{-1}) s +
      P\widetilde{a}Q c P T_r^{-1}(c)T(\widetilde{a}^{-1}) s \\
    & \phantom{=} -P\widetilde{a}QcP T_r^{-1}(c)T(\widetilde{a}^{-1}) s
    =P\widetilde{a}P\widetilde{a}^{-1} P s,
\end{align*}
 which implies the relation
  \begin{equation}\label{eq4}
I_1+I_5+I_6=P\widetilde{a}P\widetilde{a}^{-1} P s +
P\widetilde{a}Q\widetilde{a}^{-1}P s=s.
\end{equation}

Further,
\begin{align*}
I_2+I_3& = -P\widetilde{b}P JQcP T_r^{-1}(c)T(\widetilde{a}^{-1}) s
 + P\widetilde{b}P JQ\widetilde{a}^{-1}P s \\
 &= -JQbQcP T_r^{-1}(c)T(\widetilde{a}^{-1}) s +JQbQ
 \widetilde{a}^{-1} Ps \\
 &=-JQbcP T_r^{-1}(c)T(\widetilde{a}^{-1}) s +JQbQ
 \widetilde{a}^{-1} Ps\\
 & \phantom{=}\,\,+ JQb \widetilde{a}^{-1}Ps
 -JQbP\widetilde{a}^{-1}Ps\\
 &= -JQaP T_r^{-1}(c)T(\widetilde{a}^{-1}) s +JQdPs,
\end{align*}
 and consequently
\begin{equation}\label{eq5}
I_2+I_3+I_4=JQdPs.
\end{equation}
Comparing \eqref{eq4} and \eqref{eq5}, one obtains the claim for the
function $\vp_+(s)$. The case of function $\vp_-(s)$ is considered
analogously.
  \rbx

   \begin{prop}\label{c2}
 Let $(c,d)$ be the subordinated pair for a matching pair $(a,b)\in L^\infty \times L^\infty$.
 If the operator $T(c)$ is right-invertible, then
 \begin{equation}\label{direct}
\begin{aligned}
   \ker(T(a)+H(b))& =\vp_+(\im \mathbf{P}_d^+) \dotplus\im \mathbf{P}_c^-, \\
\ker(T(a)-H(b))& =\vp_-(\im \mathbf{P}_d^-) \dotplus\im
\mathbf{P}_c^+.
\end{aligned}
\end{equation}
 \end{prop}
  \textbf{Proof.}
As was mentioned in the proof of Lemma \ref{l2}, for $s\in \ker
T(d)$ one has
\begin{align*}
 (1/2)(P \widetilde{b}P+P \widetilde{a} J P)\vp^+(s)&=
\mathbf{P}_d^+ s,\\
 (1/2)(P \widetilde{b}P-P \widetilde{a} J P)\vp^-(s)&=
 \mathbf{P}_d^- s.
\end{align*}
On the other hand, if $s\in \ker T(c)$, then
 \begin{equation*}
(1/2)(\vp^+(s),\vp^-(s))^T=(\mathbf{P}_c^-s,\mathbf{P}_c^+s)^T \in
\ker\diag(T(a)+H(b),T(a)-H(b)),
 \end{equation*}
and it is easily seen that
\begin{align*}
(P \widetilde{b}P+P \widetilde{a} J P)(s -JQcPs)&=0,\\
(P \widetilde{b}P-P \widetilde{a} J P)(s +JQcPs)&=0.
\end{align*}
Using these relations one can observe that if $s\in \im
\mathbf{P}^{+}(d)$, then
 \begin{equation*}
E_2((1/2)(\vp^+(s),0)^T)=(1/2)(\vp^+(s),s)^T,
 \end{equation*}
and if $s\in \im \mathbf{P}_c^-$, then
 \begin{equation*}
E_2((1/2)(\vp^+(s),0)^T)= (1/2)(\vp^+(s),0)^T=(\mathbf{P}_c^-s,0)^T.
 \end{equation*}
Thus
 \begin{equation*}
E_2(\vp^+(\im\mathbf{P}_d^+))\cap E_2(\vp^+(\im
\mathbf{P}_c^-))=\{0\}.
 \end{equation*}
But $E_2$ is an isomorphism, hence
 \begin{equation*}
\vp^+(\im\mathbf{P}_d^+)\cap \im \mathbf{P}_c^-=\{0\},
 \end{equation*}
and it is clear that $\vp^+(\im\mathbf{P}_d^+)$ and $\im
\mathbf{P}_c^-$ are closed subspaces of $\ker (T(a)+H(b))$.
Moreover, the direct sum  of $\vp^+(\im\mathbf{P}_d^+)$ and $\im
\mathbf{P}_c^-$ is a closed subspace. In order to show that
$Y:=\vp^+(\im\mathbf{P}_d^+)\oplus \im \mathbf{P}_c^-$ is
$\ker(T(a)+H(b))$ we have to show that for $s\in \mathbf{P}_d^-$,
the element $\vp^+(s)$ belongs to $Y$. Thus assume that
$s\in\mathbf{P}_d^-$ and consider the element
 \begin{equation*}
E_2 ((\vp^+(s),0)^T)=(\vp^+(s),0)^T \in \ker T(V(a,b)).
 \end{equation*}
By Proposition \ref{p3.4}, we have $\vp^+(s)\in \ker T(c)$.
Moreover, $\vp^+(s)\in \mathbf{P}_c^-$ since otherwise
$E_1((\vp^+(s),0)^T)\notin \ker(T(a)+H(b))$, which is a
contradiction.

The related  result for  $\ker(T(a)-H(b))$ can be proved
analogously.

  \rbx

Despite the fact that all results of this section are formulated for
operators acting on $H^p$-spaces, they remain true for Toeplitz plus
Hankel operators on $l^p(\sZ_+)$ and for Wiener--Hopf plus Hankel
operators on $L^p(\sR^+)$, $1\leq p< \infty$.

\section{Coburn--Simonenko theorem for particular classes
of Toeplitz plus Hankel operators\label{s4}}

Coburn--Simonenko theorem for Toeplitz operators claims that if
$a\in L^\infty$ and is different from the zero element, then
 $$
\min(\dim\ker T(a), \dim\coker T(a))=0.
 $$
It turns out that, in general, such a statement for the operators
$T(a)+H(b)$, $a,b\in L^\infty$ is not true. Nevertheless, the
results of Section \ref{s3} allow one to single out certain classes
of Toeplitz plus Hankel operators where some kind of
Coburn--Simonenko theorem still remains in force.

 \begin{thm}\label{t2.1}
Let $a_j\in GL^\infty$, $j=1,2,3,4$ and let $A$ refer to one of the
four operators $T(a_1)-H(a_1 t^{-1})$, $T(a_2)+H(a_2 t )$,
$T(a_3)+H(a_3)$, $T(a_4)-H(a_4)$. Then $\ker A=0$ or $\coker A=0$.
 \end{thm}

  \textbf{Proof.}
Let us start with the operator $T(a_2)+H(ta_2)$. The duo
$(a_2,ta_2)$ constitutes a matching pair with the subordinated pair
$(t^{-1}, d)$, where $d=a_2 \widetilde{a}_2^{-1}t$.  Note that the
constant function $\mathbf{e}:=\mathbf{e}(t)=1$, $t\in\sT$, belongs
to both spaces $\ker T(t^{-1})$ and $\ker (T(a_2)-H(a_2 t))$.
Suppose now that $\dim \ker T(d)>0$. By Proposition \ref{c2}, we get
that
\begin{equation*}
\dim\ker(T(a_2)+H(a_2 t))=\dim \im \mathbf{P}^+_d.
\end{equation*}
Moreover, Coburn--Simonenko theorem gives  that
\begin{equation*}
 \coker T(t^{-1})=\coker T(d)=\{0\}.
 \end{equation*}
Factorization \eqref{eq3} entails that the cokernel of $T(V(a_2,a_2
t))$ is trivial. By \eqref{eqTV} the cokernel of the diagonal
operator $\diag (T(a_2)+H(a_2 t),T(a_2)-H(a_2 t))$ is also trivial,
and hence so is the cokernel of $T(a_2)+H(a_2 t)$. On the other
hand, if $\dim\ker T(d)=0$, we again can use Proposition \ref{c2} to
conclude that
 \begin{equation*}
\ker (T(a_2)+H(a_2 t))=\{0\}.
 \end{equation*}
Consider now the operator $T(a_1)-H(t^{-1}a_1)$. The representation
\eqref{eq2.2} implies that
\begin{equation}\label{eq4.1}
T(a_1)-H(a_1 t^{-1})= (T(a_1 t^{-1})-H((a_1 t^{-1})t))\cdot T(t).
\end{equation}
Setting $a_2:=a_1t^{-1}$, one rewrites the first operator in the
right-hand side of \eqref{eq4.1} as
\begin{equation*}
T(a_1 t^{-1})-H((a_1 t^{-1})t)=T(a_2)- H(a_2 t).
\end{equation*}
But the operators of the form $T(a_2)- H(a_2 t)$ have been already
considered. In particular, the element $\mathbf{e}$ belongs to the
kernel of the operator $T(a_2)- H(a_2 t)$. Now let $\ker
T(d)=\{0\}$. Since $\mathbf{e}\notin \im T(t)$, the relation
\eqref{eq4.1} implies that
\begin{equation*}
\ker(T(a_1)-H(a_1 t^{-1}))=\{0\}.
\end{equation*}
 Now let $\dim\ker T(d)>0$. Then
\begin{equation*}
\ker(T(a_2)-H((a_2 t))=\vp_-(\im \mathbf{P}_d^-)
+\Span\{\mathbf{e}\}
\end{equation*}
and
 $$
 \coker (T(a_2)-H(a_2 t))=\{0\}.
 $$
Note that the functions  $s\in \vp_-(\im \mathbf{P}_d^-) +
\Span\{\mathbf{e}\}$ can be rewritten as
\begin{equation*}
\{ (s-\widehat{s}_0)+ \widehat{s}_0:s\in \ker(T(a_2)-H(a_2t))\},
\end{equation*}
where $\widehat{s}_0$ is the zero Fourier coefficient of the
function $s$. We have $s-\widehat{s}_0 \in \im T(t)$, so if
$\widehat{s}_0\neq0$, then $\widehat{s}_0\notin \im T(t)$, and using
\eqref{eq4.1} we get
$$
\coker(T(a_1)-H(a_1 t^{-1}))=\{0\}.
$$
The remaining operators $T(a_3)+H(a_3)$ and $T(a_4)-H(a_4)$ can be
considered analogously.
  \rbx

 \begin{cor}\label{c9}
If the conditions of Theorem \ref{t2.1} are satisfied and if, in
addition, some of the corresponding operators is generalized
invertible, then it is one sided invertible.
 \end{cor}

Recall that if a linear operator is Fredholm, then it is generalized
invertible.

 \begin{rem}\label{rr3}
All of the above operators have been previously considered in
literature \cite{BE2004, BE2006, Ehr:2004h}. Nevertheless, the
proofs presented here are essentially simpler and the results of
Theorem \ref{t2.1} and Corollary \ref{c9} are more general.
\end{rem}

Of course, now one can ask how the operators $T(a_1) +H( a_1
t^{-1})$ and $T(a_2)-H(a_2 t)$ behave. For such operators, the
situation is more complicated and the assertion of Theorem
\ref{t2.1} remains valid only under additional assumptions. Consider
for instance the operator $T(a_1) +H( a_1 t^{-1})$ which was
previously studied in \cite{Ehr:2004h} in case of piecewise
continuous function $a$. Analogously to the proof of Theorem
\ref{t2.1} we set $a_2=a_1 t^{-1}$ and represent this operator as
the product of two operators,
 $$
T(a_1) +H( a_1 t^{-1})=(T(a_2) +H(a_2 t)) T(t).
 $$
We already know that $\ker T(t^{-1})$ does not affect the kernel of
the operator $T(a_2) +H(a_2 t)$. Indeed, one has $\ker T(t^{-1})=\sC
\mathbf{e}$ and $\sC \mathbf{e} \subset \ker(T(a_2) -H(a_2 t))$ and
Proposition \ref{c2} implies the claim. If $\ker T(d)=\{0\}$, then
using Proposition \ref{c2} once more we obtain that $\ker
(T(a_2)+H(a_2 t))=\{0\}$ and therefore $\ker(T(a_1) +H( a_1
t^{-1}))=\{0\}$. Assume now that $\dim\ker T(d)>0$. By Proposition
\ref{c2}, one has $\ker(T(a_2)+H(a_2 t))=\vp_+(\im \mathbf{P}_d^+)$.
Moreover, $\coker(T(a_2)+H(a_2 t))=\{0\}$ as it was shown in Theorem
\ref{t2.1}. For the sake of simplicity, suppose now that
$T(a_2)+H(a_2 t)$ is an onto operator. If $\ker(T(a_2)+H(a_2 t))
\subset \im T(t)$, then $\dim\coker (T(a_1)+H(a_1 t^{-1}))=1$. If
$\im\mathbf{P}_d^+ \neq \{0\}$, then $T(a_1)+H(a_1 t^{-1})$ has
non-trivial kernel and cokernel. If $\ker(T(a_2)+H(a_2t)$ is not
contained in $\im T(t)$, we have $\coker (T(a_1)+H(a_1
t^{-1}))=\{0\}$, the proof of which is similar to the proof of the
assertion $\coker (T(a_1)-H(a_1 t^{-1}))=\{0\}$  in the proof of
Theorem \ref{t2.1}.

It is worth noting that the study of the operator $T(a_1)+H(a_1
t^{-1})$ in \cite{Ehr:2004h} is more involved and lengthy.

 \begin{rem}\label{r3}
The results of this section are also valid for related bounded
Toeplitz plus Hankel operators considered on $l^p(\sZ_+)$. The case
of Wiener--Hopf plus Hankel integral operators requires more work
and will be published elsewhere.
 \end{rem}

\section{Kernel decomposition for a class of Toep\-litz operators\label{s5}}

In this section we establish certain properties of the kernels of
Toeplitz operators which are needed in what follows. More precisely,
we present the kernel decomposition for Fredholm Toeplitz operators
$T(g)$ with symbols $g\in L^\infty$ satisfying the relation $g
\widetilde{g}=1$. As we already know (see Corollary \ref{c1}), the
kernel of the operator $T(g)$ can be represented in the form
\begin{equation*}
\ker T(g)=\im \mathbf{P}_g^{-}\dotplus\im\mathbf{P}_g^{+}.
\end{equation*}
It turns out that the spaces $\im \mathbf{P}_g^{-}$ and $\im
\mathbf{P}_g^{+}$ possess nice bases and this fact is actively used
in the forthcoming sections. In particular, the dimensions of the
subspaces $\im \mathbf{P}_g^{+}$ and  $\im \mathbf{P}_g^{-}$ can be
determined.

If the operator $T(g): H^p\to H^p$ is Fredholm, then by Theorem
\ref{t1}, the function $g$ admits a Wiener--Hopf factorization
$g=g_- t^n g_+$ in $H^p$. Recall that $g_-(\infty)=1$.

 \begin{prop}\label{p2}
Let $g\in L^\infty$ be a function satisfying the condition $g
\widetilde{g}=1$ and such that the operator $T(g):H^p\to H^p$ is
Fredholm. Then
\begin{equation}\label{eqPM}
g_+(0)=\pm 1.
\end{equation}
  \end{prop}
 \textbf{Proof.}
Without loss of generality, we can assume that the index of $T(g)$
is equal to $0$, so the operator $T(g)$ is invertible on $H^p$. By
\cite[Proposition 7.19(c)]{BS}, the operator $T(\widetilde{g})$ is
invertible on $H^q$, $p^{-1}+q^{-1}=1$. Then the function
$\widetilde{g}=g^{-1}$ admits  Wiener--Hopf factorization in $H^q$,
 \begin{equation*}
\widetilde{g}=\widetilde{g}_+ \widetilde{g}_-=g_-^{-1} g_+^{-1}.
 \end{equation*}
It follows that
\begin{equation}\label{eqH}
\widetilde{g}_-g_+=g_-^{-1} \widetilde{g}_+^{-1}.
\end{equation}
Note that the left--hand side of \eqref{eqH} belongs to the space
$H^1(\sT)$ whereas the right--hand side is in $\overline{H^1(\sT)}$.
Taking into account the relation $H^1(\sT)\cap
\overline{H^1(\sT)}=\sC$, one obtains that there is a constant
$\xi\in\sC$ such that
 \begin{equation*}
\widetilde{g}_-g_+=g_-^{-1} \widetilde{g}_+^{-1}=\xi,
 \end{equation*}
so
\begin{equation*}
g_+=\xi \widetilde{g}_-^{-1}, \quad g_-^{-1}=\xi \widetilde{g}_+.
\end{equation*}
 Moreover, recalling the identities $g_-(\infty)=g_-^{-1}(\infty)=1$ and
 $\widetilde{g}_-^{-1}(0)=g_-(\infty)=1$ and
 $\widetilde{g}_+(\infty)=g_+(0)$ we get $g_+(0)=\xi$ and
 $\xi^2=1$, which implies \eqref{eqPM}.
 \rbx

 \begin{defn}\label{defn1}
Let $g\in L^\infty$ satisfy the condition $g \widetilde{g}=1$ and
admit the Wie\-ner--Hopf factorization $g=g_-t^n g_+$,
$g_-(\infty)=1$ in $H^p$. The value $g_+(0)$ is called the
factorization signature of the function $g$ and is denoted by
$\boldsymbol\sigma(g)$.
 \end{defn}

 \begin{cor}\label{c3}
If $g\in L^\infty$ is a matching function such that the operator
$T(g)$ is Fredholm, then the factors $g_-$ and $g_+$ in the related
Wiener--Hopf factorization of $g$ satisfy the relations
 \begin{align*}
g_+=\boldsymbol\sigma (g)\widetilde{g}_-^{-1}, \quad
g_-=\boldsymbol\sigma (g)\widetilde{g}_+^{-1}.
\end{align*}
 \end{cor}
Indeed, these relations immediately follow from the proof of
Proposition \ref{p2}.

The following result plays a crucial role in this paper.

\begin{thm}\label{t2}
Let $g\in L^\infty$ be a matching function such that the operator
$T(g):H^p\to H^p$ is Fredholm and $n:=\ind T(g)>0$. If $g=g_-t^{-n}
g_+$, $g_-(\infty)=1$ is the corresponding Wiener--Hopf
factorization of $g$ in $H^p$, then  the following systems of
functions $\cB_{\pm}(g)$ form bases in the spaces $\im
\mathbf{P}_g^{\pm}$:
  \begin{enumerate}
    \item If $n=2m$, $m\in \sN$, then
    \begin{equation*}
\cB_{\pm}(g):= \{g_+^{-1}(t^{m-k-1}\pm \boldsymbol\sigma
(g)t^{m+k}):k=0,1,\cdots, m-1
 \},
    \end{equation*}
     and
     \begin{equation*}
      \dim\im \mathbf{P}_g^{\pm}=m.
     \end{equation*}
    \item If $n=2m+1$, $m\in \sZ_+$, then
    \begin{equation*}
\cB_{\pm}(g):= \{g_+^{-1}(t^{m+k}\pm \boldsymbol\sigma
(g)t^{m-k}):k=0,1,\cdots, m-1
 \}\setminus \{0\},
    \end{equation*}
    \begin{equation*}
     \dim\im \mathbf{P}_g^{\pm}=m +\frac{1\pm\boldsymbol\sigma(g)}{2}
     \,,
    \end{equation*}
and the zero element belongs only to one of the sets $\cB_{+}(g)$ or
$\cB_{-}(g)$. Namely, for $k=0$ one of the terms
$t^m(1\pm\boldsymbol\sigma(g))$ is equal to zero.
   \end{enumerate}
\end{thm}
 \textbf{Proof.}
 It is easily seen that the restriction of the operators $Pg_-I$
 and $Pg_-^{-1}I$ on $\ker T(t^{-n})=\Span \{\mathbf{e}, \cdots,
 t^{n-1} \}$ map $\ker T(t^{-n})$ into $\ker T(t^{-n})$ and on the
 space $\ker T(t^{-n})$ the above operators are inverses to each
 other.

Clearly, the elements $s_j=Pg_-t^j$, $j=0,1,\cdots, n-1$ are again
in $\ker T(t^{-n})$ and
\begin{equation*}
T^{-1}(g_0)s_j=g_+^{-1} Pg_-^{-1}s_j=g_+^{-1}t^j,
\end{equation*}
where $g_0:=gt^n$.

Note that $T^{-1}(g_0)=g_+^{-1}Pg_-^{-1}$ and $T^{-1}(g_0)s_j\in
\ker T(g)$. Moreover, the set $\{T^{-1}(g_0)s_j:j=0,\cdots,n-1 \}$
constitutes a basis in $\ker T(g)$. Now one can consider the
expression
 \begin{align*}
JQgPT^{-1}(g_0)s_j &= JQg_+ t^{-n}g_- g_+^{-1} Pg_-^{-1} s_j \\
    &= JQt^{-n}g_-t^j = Pt^n\widetilde{g}_- t^{-j-1}= P
    t^{n-j-1}\widetilde{g}_-.
\end{align*}
Corollary \ref{c3} shows that $\widetilde{g}_-=\boldsymbol\sigma (g)
g_+^{-1}$, which leads to the relation
\begin{equation*}
\mathbf{P}_g^{\pm} T^{-1}(g_0) s_j= \frac{1}{2}g_+^{-1}(t^j\pm
\boldsymbol\sigma(g)t^{n-j-1}), \quad j=0,1,,\cdots,n-1.
\end{equation*}

Now let $n=2m$, $m\in\sN$. If $j\in\{0,1,\cdots,m-1\}$, then this
$j$ can be written as $j=m-k-1$ with some $k\in\{0,1,\cdots,m-1\}$
and vice versa. Hence
\begin{equation*}
t^j\pm \boldsymbol\sigma(g)t^{n-j-1} = t^{m-k-1}\pm
\boldsymbol\sigma(g)t^{m+k}, \quad j=0,1,\cdots,m-1.
 \end{equation*}
On the other hand, if $j\geq m$, then $j$ can be rewritten as
$j=m+k$ for a $k\in\{0,1,\cdots, m-1\}$, and $t^j\pm
\boldsymbol\sigma(g)t^{n-j-1} = t^{m+k}\pm
\boldsymbol\sigma(g)t^{m-k-1}$. Thus one obtains, maybe up to the
factor $-1$,  the same function system that was found for
$j\in\{0,1,\cdots, m-1\}$. So we conclude that if $n=2m$, then
 \begin{equation*}
\dim \im \mathbf{P}_g^{\pm}=m,
 \end{equation*}
and assertion (i) is shown.  Assertion (ii) can be proved similarly.
 \rbx

Let us emphasize that, in general, the determination of the
factorization signature $\boldsymbol\sigma(g)$, $g\in L^\infty$ is a
difficult problem. Nevertheless,  we can provide sufficient
conditions which allow one to obtain $\boldsymbol\sigma(g)$. Note
that if $g\in L^\infty$ is continuous at the points $\pm 1\in\sT$,
then the function $g$ can take only two values $-1$ or $1$ at the
points mentioned, i.e. $g(-1),g(1)\in \{-1,1\}$.

  \begin{prop}\label{p3}
Let $g\in L^\infty$ be a matching function such that
\begin{enumerate}
    \item The operator $T(g)$ is invertible on $H^p$.
    \item The function $g$ is continuous at the point $1$ or $-1$.
\end{enumerate}
Then  $\boldsymbol\sigma(g)=g(1)$ or  $\boldsymbol\sigma(g)=g(-1)$,
respectively.
  \end{prop}
 \textbf{Proof.} Assume for definiteness that the function $g$ is continuous at the
point $1$. Then it admits a Wiener--Hopf factorization in $H^p$,
\begin{equation*}
g=g_- g_+, \quad g(\infty)=1,
\end{equation*}
and $g(1)=\pm 1$. By Corollary \ref{c3}  one has
\begin{equation*}
g=\boldsymbol\sigma(g)\widetilde{g}_+^{-1} g_+ \, .
\end{equation*}
Approximate the function $g$ as follows. For a given $\ve>0$ chose
an arc of $\sT$ with endpoints $e^{i\xi_0}$ and $e^{-i\xi_0}$ such
that the point $1$ belongs to this arc and such that the function
\begin{equation*}
g_{\ve}(t)=\left \{
 \begin{array}{cc}
   g(1) & \text{ if }\, t= e^{i\theta}, \, \theta\in (-\xi_0,\xi_0) \\
   g(t) & \text{ otherwise }\\
 \end{array}
 \right . ,
\end{equation*}
satisfies the condition
\begin{equation*}
||g-g_{\ve}||<\ve .
\end{equation*}
If $\ve$ is small enough, then the operator $T(g_{\ve})$ is also
invertible and
 \begin{equation}\label{AF}
g_{\ve}=\boldsymbol\sigma(g_{\ve})(\widetilde{g}_{\ve})_+^{-1}
(g_{\ve})_+ .
\end{equation}
Since $g_{\ve}$ is H\"older continuous in a neighbourhood of the
point $1$, \cite[Corollary 5.15]{LiSp:1987} shows that the functions
$(g_{\ve})_+$ and $(\widetilde{g}_{\ve})_+^{-1}$ are also H\"older
continuous in a neighbourhood of the point $1\in \sT$. From relation
\eqref{AF} one then obtains
\begin{equation*}
g_{\ve}(1)=\boldsymbol\sigma(g_{\ve}).
\end{equation*}
Additionally assume that $\ve$ is so small that
  \begin{equation}\label{Eq1}
  ||T^{-1}(g)- T^{-1}(g_{\ve})|| <1.
 \end{equation}
The equations $T(g)h=1$ and $T(g_{\ve}k=1$ are uniquely solvable and
using Wiener--Hopf factorizations
 \begin{align*}
    g=g_- \, g_+, &\quad g_-(\infty)=1,\\
 g_{\ve}=(g_{\ve})_- \, g(_{\ve})_+,& \quad (g_{\ve})_-(\infty)=1,
\end{align*}
one obtains
\begin{align*}
    h&=T^{-1}(g) \mathbf{e}= g_+^{-1}P g_-^{-1}\mathbf{e}= g_+^{-1}\\
 k&=T^{-1}(g_{\ve}) \mathbf{e}= (g_{\ve})_+^{-1}P (g_{\ve})_-^{-1}\mathbf{e}= (g_{\ve})_+^{-1},
\end{align*}
where $\mathbf{e}$ is the constant function
$\mathbf{e}:=\mathbf{e}(t)=1$, $t\in\sT$. Now the inequality
\eqref{Eq1} implies the estimate
 \begin{equation*}
||g_+^{-1}-(g_{\ve})_+^{-1}||_p<\ve .
\end{equation*}
Consequently, equations $g_+^{-1}(0)=\pm1$,
$(g_{\ve})_+^{-1}(0)=\pm1$ lead to the relation
$g_+^{-1}(0)=(g_{\ve})_+^{-1}$, i.e.
$\boldsymbol\sigma(g)=\boldsymbol\sigma(g_{\ve})$. Hence,
\begin{equation*}
\boldsymbol\sigma(g)=\boldsymbol\sigma(g_{\ve})=g_{\ve}(1)=g(1),
\end{equation*}
which completes the proof.
  \rbx

\begin{cor}\label{c4}
Let $g\in L^\infty$  be a matching function satisfying condition
(ii) of Proposition \ref{p3}. If the operator $T(g):H^p\to H^p$ is
Fredholm and $n:=\ind T(g)$, then $\boldsymbol\sigma(g)=g(1)$ if $g$
is continuous at the point $1$ and $\boldsymbol\sigma(g)=g(1)$.
    \item If $n$ is even, then   $\boldsymbol\sigma(g)=g(-1)(-1)^n$ if
    $g$ is continuous at the point $-1$.
 \end{cor}

 \section{Toeplitz plus Hankel operators with Fredholm matching pair\label{s6}}

In this section the structure of the kernel and cokernel of Toeplitz
plus Hankel operator $T(a)+H(b)$ is described. The operators in
question are studied under the condition that their generating
functions $a, b\in L^\infty$ constitute a Fredholm matching pair
$(a,b)$. Recall that if a matching pair $(a,b)$ is Fredholm, then it
follows from \eqref{eqTV} and \eqref{eq3} that $T(a)+H(b)$ and
$T(a)-H(b)$ are Fredholm operators. Set $\kappa_1:=\ind T(c)$,
$\kappa_2:=\ind T(d)$ and let $\sZ_-$ refer to the set of all
negative integers.
 \begin{thm}\label{t3}
Assume that  $(a,b)\in L^\infty\times L^\infty$ is a Fredholm
matching pair. Then
 \begin{enumerate}
    \item If $(\kappa_1,\kappa_2)\in \sZ_+\times \sZ_+$, then the
    operators $T(a)+H(b)$ and $T(a)-H(b)$ are invertible from the
    right and
    \begin{align*}
\ker (T(a)+H(b))&=\im\mathbf{P}_c^{-}\dotplus \vp_+(\im \mathbf{P}_d^{+}), \\
\ker (T(a)-H(b))&=\im\mathbf{P}_c^{+}\dotplus \vp_-(\im
\mathbf{P}_d^{-}),
\end{align*}
where the spaces $\im\mathbf{P}_c^{\pm}$ and $\im\mathbf{P}_d^{\pm}$
are described in Theorem \ref{t2}, and the mappings $\vp_{\pm}$ are
defined by \eqref{eqphi}.

     \item  If $(\kappa_1,\kappa_2)\in (\sZ\setminus \sN)\times (\sZ\setminus \sN)$, then the
    operators $T(a)+H(b)$ and $T(a)-H(b)$ are invertible from the
    left and
    \begin{align*}
\coker (T(a)+H(b))&=\im
\mathbf{P}_{\overline{d}}^{-}\dotplus \vp_+(\im\mathbf{P}_{\overline{c}}^{+}), \\
\coker (T(a)-H(b))&=\im \mathbf{P}_{\overline{d}}^{+}\dotplus
\vp_-(\im\mathbf{P}_{\overline{c}}^{-}),
\end{align*}
and $\im \mathbf{P}_{\overline{d}}^{\pm}=\{0\}$ for $\kappa_2=0$.

   \item If $(\kappa_1,\kappa_2)\in \sZ_+\times \sZ_-$, then
    \begin{align*}
\ker (T(a)+H(b))=\im\mathbf{P}_c^{-},  & \quad  \coker (T(a)+H(b))=\im \mathbf{P}_{\overline{d}}^{+}, \\
\ker (T(a)-H(b))=\im\mathbf{P}_c^{+}  & \quad  \coker
(T(a)-H(b))=\im \mathbf{P}_{\overline{d}}^{-}.
\end{align*}

 \end{enumerate}
 \end{thm}

 \textbf{Proof.}
Let us note that all results concerning the kernels of the operators
mentioned follow from Theorem \ref{t2}. Considering the cokernels of
the corresponding operators, we recall that $\coker (T(a)\pm
H(b)):=\ker (T(a)\pm H(b))^*$. Moreover, $(T(a)\pm
H(b))^*=T(\overline{a})\pm H(\widetilde{\overline{b}})$ and
$(\overline{a}, \widetilde{\overline{b}})$ is again a matching pair
with the subordinated pair $(\overline{d}, \overline{c})$. Further,
if $c=c_-t^{-\kappa_1} c_+$ is the Wiener--Hopf factorization of $c$
in $H^p$, then $\overline{c}=(\boldsymbol\sigma(c)
\overline{c}_+)t^{\kappa_1}(\boldsymbol\sigma(c)\overline{c}_-)$ is
the related Wiener--Hopf factorization of $\overline{c}$ in $H^q$
and $\overline{c}_-\in H^p$, $\overline{c_-^{-1}}\in H^q$,
$\overline{c}_+\in \overline{H^q}$, $\overline{c_-^{-1}}\in
\overline{H^p}$, $p^{-1}+q^{-1}=1$, and
$\boldsymbol\sigma(\overline{c})=\boldsymbol\sigma(c)$. Of course,
since the function $\overline{d}$ admits a similar factorization,
cokernel description can be obtained directly from the previous
results for the kernels of Toeplitz plus Hankel operators.
 \rbx

It remains to consider the case  $(\kappa_1,\kappa_2)\in \sZ_-\times
\sZ_+$. This situation is more involved and factorization
\eqref{eq3} already indicates that for $\kappa_2>0$, the kernel
dimension of $\diag(T(a)+H(b),T(a)-H(b))$ may be smaller than
$\kappa_2$. In order to prepare our next theorem, choose an
$n\in\sN$ such that
\begin{equation*}
1\geq 2n+\kappa_1\geq 0.
\end{equation*}
Such an $n$ is uniquely defined and
\begin{equation*}
2n+\kappa_1 =\left\{%
\begin{array}{ll}
   0, & \hbox{if\;} \kappa_1 \; \hbox{is even,}  \\
    1, &\hbox{if\;} \kappa_1 \; \hbox{is odd.} \\
\end{array}%
\right.
\end{equation*}
Now the operators $T(a)\pm H(b)$ can be represented in the form
 \begin{equation}\label{eq6.1}
T(a)\pm H(b)= ( T(at^{-n})\pm H(bt^n))T(t^n).
\end{equation}
Note that $(at^{-n}, bt^{n})$ is a matching pair with the
subordinated pair $(ct^{-2n}, d)$. Therefore, the operators
$T(at^{-n})\pm H(bt^n)$ are subject to assertion (i) of Theorem
\ref{t3}. Thus they are right-invertible, and if $\kappa_1$ is even,
then
 \begin{equation}\label{eq6.2}
\begin{aligned}
    \ker (T(at^{-n})+ H(bt^n) )=\vp_+(\im \mathbf{P}_d^{+}),\\
 \ker (T(at^{-n})- H(bt^n) )=\vp_-(\im \mathbf{P}_d^{-}),
\end{aligned}
\end{equation}
and if $\kappa_1$ is odd, then
\begin{equation}\label{eq6.3}
\begin{aligned}
    \ker (T(at^{-n})+ H(bt^n) )= \frac{1-\boldsymbol\sigma(c)}{2}c_+^{-1}\sC
\dotplus  \vp_+(\im \mathbf{P}_d^{+}),\\
 \ker (T(at^{-n})- H(bt^n) )=\frac{1+\boldsymbol\sigma(c)}{2}c_+^{-1}\sC
 \dotplus \vp_-(\im \mathbf{P}_d^{-}),
\end{aligned}
\end{equation}
where the mappings $\vp_{\pm}$ depend on the functions $at^{-n}$ and
$bt^n$.

 \begin{thm}\label{t4}
Let $(\kappa_1,\kappa_2)\in \sZ_-\times \sZ_+$. Then
\begin{enumerate}
    \item If $\kappa_1$ is odd, then
   \begin{multline*}
\hspace*{-6mm}\ker(T(a) \pm  H(b))\!  = \!T(t^{-n})\left (\left\{
 \frac{1\mp\boldsymbol\sigma(c)}{2}c_+^{-1}\sC
 \dotplus \vp_{\pm}(\im \mathbf{P}_d^{\pm})\right\} \cap \im T(t^n)\right )\\
 =\left\{\psi\in \{ T(t^{-n})u\}:   u\in  \left
\{\frac{1\mp\boldsymbol\sigma(c)}{2}c_+^{-1}\sC
 \dotplus \vp_{\pm}(\im \mathbf{P}_d^{\pm})\right \} \right . \\
 \left. \phantom{\frac{1\mp\boldsymbol\sigma(c)}{2}}  \text{and }\,
   \widehat{u}_0=\cdots =\widehat{u}_{n-1}=0,   \right \},
\end{multline*}
where $\widehat{u}_k, k=0,1,\cdots, n-1$ are the Fourier
coefficients of the function $u$, and the mappings $\vp_{\pm}$
depend on the functions $at^{-n}$ and $bt^n$.
    \item If $\kappa_1$ is even, then
     \begin{align*}
 &\ker(T(a)\pm  H(b)) =  T(t^{-n})\left ( \left\{\vp_{\pm}(\im \mathbf{P}_d^{\pm})\right\}
 \cap \im T(t^n) \right )\\
 &= \left\{ \psi\in \{ T(t^{-n})u\}:
u\in  \vp_{\pm}(\im \mathbf{P}_d^{\pm}) \, \text{and }\,
 \widehat{u}_0=\cdots =\widehat{u}_{n-1}=0\right \},
 \end{align*}
and the mappings $\vp_{\pm}$ again depend on $at^{-n}$ and $bt^n$.
\end{enumerate}
\end{thm}

 \textbf{Proof.}
It follows immediately from representations
\eqref{eq6.1}--\eqref{eq6.3}. \rbx

Theorem \ref{t4} can also be used to derive representations of the
cokernel of the operator $T(a)\pm H(b)$ in the situation where
$(\kappa_1,\kappa_2)\in \sZ_-\times \sZ_+$. Indeed, recalling that
$(T(a)\pm H(b))^*=T(\overline{a})\pm H(\widetilde{\overline{b}})$,
and $(\overline{d}, \overline{c})$ is the subordinated pair for
$(\overline{a},\widetilde{\overline{b}})$, one can note that the
operators $T(\overline{d})$ and $T(\overline{c})$ are also Fredholm
and
\begin{equation*}
\ind T(\overline{d})=-\kappa_2, \quad \ind
T(\overline{c})=-\kappa_1,
\end{equation*}
so $(-\kappa_2,\kappa_1)\in \sZ_-\times \sZ_+$. Therefore, Theorem
\ref{t4} applies and we can formulate the following result.

 \begin{thm}\label{t5}
Let $(\kappa_1,\kappa_2)\in \sZ_-\times \sZ_+$, and let $m\in\sN$
satisfy the requirement
\begin{equation*}
1\geq 2m-\kappa_2\geq0.
\end{equation*}
Then
\begin{enumerate}
    \item If $\kappa_2$ is odd, then
     \begin{equation*}
\hspace*{-6.5mm}\!\!\! \coker(T(a)\pm H(b))\! = \! T(t^{-m})\!\left
(\!\left\{
 \frac{1\mp\boldsymbol\sigma(\overline{d})}{2}\overline{d_-^{-1}}\sC
 \dotplus \vp_{\pm}(\im \mathbf{P}_{\overline{c}}^{\pm})\right\} \cap \im T(t^m) \!\right ).
\end{equation*}

    \item If $\kappa_2$ is even, then
     \begin{align*}
 \ker(T(a)\pm H(b)) &= T(t^{-m})\left (\left\{\vp_{\pm}(\im \mathbf{P}_{\overline{c}}^{\pm})\right\} \cap \im
 T(t^m) \right ),
 \end{align*}
and the mappings $\vp_{\pm}$ depend on $\overline{a}t^{-m}$ and
$\widetilde{\overline{b}}t^m$.
\end{enumerate}
\end{thm}

In some cases the above approach allows one to drop the condition of
Fredholmness of the operator $T(d)$.  We are not going to pursue
this analysis here but  rather restrict ourselves to a few special
cases generalizing results of Section~\ref{s4}.

  \begin{cor}\label{c5}
Let $(a,b)\in L^\infty \times L^\infty$ be a matching pair with the
subordinated pair $(c,d)$, and let $T(c)$ be a Fredholm operator.
Then:
\begin{enumerate}
    \item If $\ind T(c)=1$ and $\boldsymbol\sigma(c)=1$, then $\ker(T(a)+H(b))=\{0\}$
     or $\coker(T(a)+H(b))=\{0\}$.
    \item If $\ind T(c)=-1$ and $\boldsymbol\sigma(c)=1$, then $\ker(T(a)-H(b))=\{0\}$ or
    $\coker(T(a)-H(b))=\{0\}$.
    \item If $\ind T(c)=0$, then $\ker(T(a)\pm H(b))=\{0\}$ or $\coker(T(a)\pm
    H(b))=\{0\}$.
\end{enumerate}
 \end{cor}
 \textbf{Proof.}
The proof of the last theorem is similar to the proof of
Theorem~\ref{t2}, but in the proofs of assertions (i) and (ii) one,
respectively, has to use the fact that $c_+^{-1}\in \ker
(T(a)+H(b))$ and $c_+^{-1}\in \ker (T(a t^{-1})-H(bt))$. These
inclusions can be verified by straightforward computations. Thus
considering, for example, the expression $(T(a
t^{-1})-H(bt))c_+^{-1}$, one obtains
  \begin{align*}
(T(a t^{-1})-H(bt))c_+^{-1}& = Pbc t^{-1}c_+^{-1}-PbtQJc_+^{-1}\\
 &=Pbc_+tc_- t^{-1}c_+^{-1}-Pbt\widetilde{c}_+^{-1} t^{-1}=0.
  \end{align*}
Note that we have used the relations $a=bc$ and
$c_-=\widetilde{c}_+^{-1}$.
 \rbx

 \begin{cor}\label{c6}
Let $b\in L^\infty$ be a matching function. If $T(\widetilde{b})$ is
a Fredholm operator, then:
\begin{enumerate}
    \item If \/ $\ind T(\widetilde{b})=1$,
    and $\boldsymbol\sigma(\widetilde{b})=1$, then $\ker(I+H(b))$  or
    $\coker(I+H(b))$ is trivial.
    \item If\/ $\ind T(\widetilde{b})=-1$,
    and $\boldsymbol\sigma(\widetilde{b})=1$, then $\ker(I-H(b))$ or $\coker(I-H(b))$ is trivial.
    \item If \/ $\ind T(\widetilde{b})=0$,
    then $\ker(I\pm H(b))$ or $\coker(I\pm H(b))$ is trivial.
\end{enumerate}
 \end{cor}
This results is a direct consequence of  Corollary \ref{c5}, since
if $b$ is a matching function, then $(1,b)$ is a matching pair with
the subordinate pair $(\widetilde{b},b)$.

 \section{$PC$-generating functions\label{s7}}

The results of the previous section can be improved if there is more
information available about the generating functions $a$ and $b$.
Thus if $a$ and $b$ are piecewise continuous functions, then
one-sided invertibility of Toeplitz plus Hankel operators on the
space $H^p$, $1<p<\infty$ can be studied in more detail. Recall that
a function $a\in L^\infty$ is called piecewise continuous if for
every $t\in \sT$ the one-sided limits $a(t+0)$ and $a(t-0)$ exist.
The set of all piecewise continuous functions is denoted by
$PC(\sT)$ or simply by $PC$. It is well-known that $PC$ is a closed
subalgebra of $L^\infty$, and any piecewise continuous function has
at most countable set of jumps. Moreover, for each $\delta>0$ the
set $S:=\{ t\in\sT:|a(t+0)-a(t-0)|>\delta \}$ is finite.

Further, let us introduce the functions
 \begin{align*}
 \nu_p(y)&:=\frac{1}{2}\left ( 1+\coth\left (\pi
\left(y+\frac{i}{p}\right ) \right )\right
    ), \quad 
 h_p(y):= \sinh^{-1}\left (\pi \left(y+\frac{i}{p}\right )
\right ),
\end{align*}
where $y\in \overline{\sR}$, and $\overline{\sR}$ refers to the
two-point compactification of $\sR$. Note that for  given points
$u,w\in \sC$, $u\neq w$ the set $\cA_p(u,w):=\{z\in\sC: z=w
\nu_p(y)+u(1-\nu_p(y)), y\in \overline{\sR}\}$ forms a circular arc
which starts at $u$ and ends at $w$ as $y$ runs through
$\overline{\sR}$. This arc $\cA_p(u,w)$ has the property that from
any point of the arc, the line segment $[u,w]$ is seen at the angle
$2\pi/(\max\{p,q\})$, $1/p+1/q=1$. Moreover, if $2<p<\infty$
($1<p<2$) the arc $\cA_p(u,w)$ is located on the right-hand side
(left-hand side) of the straight line passing through the points $u$
and $w$ and directed from $u$ to $w$. If $p=2$, the set $\cA_p(u,w)$
coincides with the line segment $[u,w]$.

Let $\tp:=\{t\in \sT:\Im t>0\}$.

 \begin{thm}\label{tTH}
 If $a,b\in PC$, then
 the operator $T(a)+H(b)$ is Fredholm if and only if the matrix
  \begin{align*}
    &\smb(T(a)+H(b))(t,y):= \nn\\[1ex]
  &  \left(%
\begin{array}{c@{\hspace{-1mm}}c}
a(t+0)\nu_p(y)+a(t-0)(1-\nu_p(y))   & \D \frac{b(t+0)-b(t-0)}{2i} \,h_p(y)\\
\D \frac{b(\overline{t}-0)-b(\overline{t}+0)}{2i} \,h_p(y)   &
 a(\overline{t}+0)\nu_p(y)+a(\overline{t}-0)(1-\nu_p(y))\\
   \end{array}%
\right) 
\end{align*}
is invertible for every $(t,y)\in \tp\times \overline{\sR}$ and the
function
 \begin{align*}
\smb (T(a)+H(b))(t,y):=& \, a(t+0)\nu_p(y)+a(t-0)(1-\nu_p(y))\\
 & +
t \, \frac{b(t+0)-b(t-0)}{2} \,h_p(y)
\end{align*}
does not vanish on $\{-1,1\}\times \overline{\sR}$.
 \end{thm}

 \begin{rem}\label{r1}
For the first time this theorem appeared in \cite{RS1990} (see also
\cite{RSS2011}). An index formula can also be established similarly
to \cite{RS:2012}. If the matching condition is satisfied, an index
formula can be derived from Theorem \ref{t6} below.
 \end{rem}

If $(a,b)\in PC\times PC$ is a matching pair, then the subordinated
pair $(c,d)\in PC\times PC$, and both operators $T(a)+H(b)$ and
$T(a)-H(b)$ are simultaneously Fredholm if and only if the matching
pair $(a,b)$ is Fredholm. This follows from the fact that
semi-Fredholm Toeplitz operators with piecewise continuous
generating functions are indeed Fredholm operators. This assertion
can be obtained immediately from \cite[Proposition 3.1]{GK1992a} and
\cite[Lemma 1]{Mu2003} and using \eqref{eqTV} and \eqref{eq3}. Let
us mention another results we need.

Let $A$ be an operator defined on all spaces $L^p$ for $1<p<\infty$.
Consider the set $A_{F}:=\{p\in(1,\infty) \; \text{such that the
operator} \; A:H^p\to H^p \; \text{is Fredholm}\}$.

 \begin{prop}[See {\cite{Sh1974}}]\label{p4}
The set $A_F$ is open. Moreover, for each connected component
$\gamma\in A_F$, the index of the operator $A:L^p\to L^p$, $p\in
\gamma$ is constant.
 \end{prop}
 For Toeplitz operators the structure of the set $A_F$ can be characterized as follows.
 \begin{prop}[See {\cite{Sp:1976}}]\label{p5}
 Let $G$ be an invertible matrix--functions with entries from $PC$,
 and let $A:=T(G)$. Then there is an at most countable subset $S_A\subset
 (1,\infty)$ with the only possible accumulation points $t=1$ and
 $t=\infty$ such that $A_F=(1,\infty)\setminus S_A$.
\end{prop}

This result can be used to describe the corresponding set $A_F$ for
Toeplitz plus Hankel operators.

 \begin{cor}\label{c7}
Let $a,b\in PC$, and let $A:=\diag(T(a)+H(b),T(a)-H(b)): H^p\times
H^p \to H^p\times H^p$. Then there is at most countable subset
$S_A\subset (1,\infty)$ with the only possible accumulation points
$t=1$ and $t=\infty$ such that $A_F=(1,\infty)\setminus S_A$.
 \end{cor}
 \textbf{Proof.}
It follows directly from Proposition \ref{p5} since
 $\diag(T(a)+H(b),T(a)-H(b))$ is Fredholm if and only if so is the
 operator $T(V(a,b))$.
\rbx

Thus if $a,b\in PC$ and the operator $T(a)+H(b)$ is Fredholm on
$H^p$, then there is an interval  $(p',p'')$ containing $p$ such
that $T(a)+H(b)$ is Fredholm on all spaces $H^r$, $r\in (p',p'')$
and the index of this operator does not depend on $r$. Moreover,
there is an interval $(p,p_0)\subset (p',p'')$, $p<p_0$ such that
$T(a)-H(b)$ is Fredholm on $H^r$, $r\in (p,p_0)$ and its index does
not depend on $r$. Now we can formulate the following result.

\begin{prop}\label{p6}
Let $T(a)+H(b):H^p\to H^p$ be a Fredholm operator and let $r\in (p,
p_0)$. Then the kernel and cokernel of the operator
$T(a)+H(b):H^r\to H^r$ coincide with the kernel and cokernel of the
same operator acting on the space $H^p$.
 \end{prop}

 \textbf{Proof.}
Let us first recall a result from  \cite{GF1974}. Assume that
$X_1,X_2$ are Banach spaces such that $X_1$ is continuously and
densely embedded into $X_2$ and $A$ is a linear bounded operator
both on $X_1$ and $X_2$. If $A$ is a Fredholm operator on each space
$X_1$ and $X_2$, and
$$
\ind A\left |_{X_1\to X_1} \right .=\ind A\left |_{X_2\to X_2},
\right .
$$
then
 \begin{align*}
\ker A\left |_{X_1\to X_1} \right .&=\ker A\left |_{X_2\to X_2}, \right .\\
\coker A\left |_{X_1\to X_1} \right .&=\coker A\left |_{X_2\to X_2},
\right .
\end{align*}
which implies the assertion.
 \rbx

Let us now formulate a result concerning the kernels and cokernels
of Hankel plus Toeplitz operators.

 \begin{thm}\label{t6}
Let $a,b\in PC$ and $(a,b)$ be a matching pair. If the operator
$T(a)+H(b):H^p\to H^p$ is Fredholm, then there is an interval
$(p,p_0)$, $p<p_0$ such that for all $r\in(p,p_0)$ the pair $(a,b)$
and both operators $T(a)\pm H(b):H^r\to H^r$ are Fredholm,
\begin{equation*}
 \begin{aligned}
\ker (T(a)+H(b))\left |_{H^r\to H^r} \right .&
 =\ker( T(a)+H(b))\left |_{H^p\to H^p}\,, \right .\\
\coker (T(a)+H(b))\left |_{H^r\to H^r} \right .&
 =\coker (T(a)+H(b))\left |_{H^p\to H^p}\, , \right .\\
\end{aligned}
\end{equation*}
and the kernel and cokernel of the operator $T(a)+H(b):H^r\to H^r$
are described by Theorems \ref{t3}--\ref{t5}.
\end{thm}
 \textbf{Proof.}
It follows immediately from Proposition \ref{p6} and previous
considerations. \rbx

 \begin{rem}
Let $a,b\in L^\infty$ be functions which possess the property fixed
in Corollary \ref{c7} for $PC$-functions. Then the results of
Theorem \ref{t6} remain true.
 \end{rem}

\section{A few remarks on factorization signature \label{s8}}

In previous sections the factorization signature has been used to
describe the kernels of Toeplitz plus Hankel operators. Therefore,
the determination of this characteristic is an important problem.
Let us consider this problem for piecewise continuous matching
functions $c$ such that the operator $T(c)$ is Fredholm. For the
sake of definiteness we assume that if $z\in \sC\setminus \{0\}$ is
a complex number, then its argument $\arg z$ is always chosen to be
in the interval $(-\pi, \pi]$. Following \cite[Section 5.35]{BS},
for $\beta\in\sC$ and $\tau\in \sT$ consider a function $\vp_{\beta,
\tau}(t)\in PC$ defined by
\begin{equation*}
\vp_{\beta, \tau}(t):=\exp\{i\beta \arg (-t/\tau) \}, \quad t\in\sT.
\end{equation*}
It is easily seen that $\vp_{\beta, \tau}$ has at most one
discontinuity, namely, a jump at the point $\tau$ and
\begin{equation*}
\vp_{\beta, \tau}(\tau+0)=e^{-i\pi\beta}, \quad \vp_{\beta,
\tau}(\tau-0)=e^{i\pi\beta}.
\end{equation*}
From \cite[Sections 5.35 and 5.36]{BS} one obtains that:
 \begin{enumerate}
    \item The matching function $c$ can be represented in one of the
following form
 \begin{equation}\label{eq7.2}
c=\vp_{\gamma_+} c_1, \, \text{ or }\, c=\vp_{\gamma_-} c_1
\end{equation}
where the function $c_1,c_2\in PC$ are continuous at the points
$t=+1$ and $t=-1$, respectively; $\vp_{\gamma_+}:=\vp_{\gamma_+,1}$,
$\vp_{\gamma_-}:=\vp_{\gamma_-,-1}$, and the parameters
$\gamma_+,\gamma_-$ are defined according to \cite[Section 5.36]{BS}
with $\re \gamma_+,\re \gamma_-\in(-1/q,1/p)$ .
    \item The operators $T(\vp_{\gamma_+})$ and $T(\vp_{\gamma_-})$
    are invertible and the operators $T(c_1)$ and $T(c_2)$ are Fredholm.
 \end{enumerate}
It is not hard to see that
 $$
\vp_{\gamma_{\pm}} \widetilde{\vp}_{\gamma_{\pm}}=1, \, \text{ and
}\,\vp_{\gamma_+}(-1)=1, \quad \vp_{\gamma_-}(1)=1.
 $$
Proposition \ref{p3} then ensures that
 $$
\boldsymbol\sigma(\vp_{\gamma_{\pm}})=1.
 $$
Moreover, $c_1$ and $c_2$ are also matching functions and by
Proposition \ref{p3} and Corollary \ref{c4} we have
 $$
\boldsymbol\sigma(c_1)=c_1(1), \quad
\boldsymbol\sigma(c_2)=c_2(-1)(-1)^n,
 $$
where $n=\ind T(c_2)=\ind T(c)$.

 \begin{thm}\label{t7}
If $c\in PC$ is a matching function which admits representation
\eqref{eq7.2}, then
 $$
\boldsymbol\sigma(c)=\boldsymbol\sigma(c_1)=\boldsymbol\sigma(c_2).
 $$
 \end{thm}

 \textbf{Proof.} Assume for definiteness that $c=\vp_{\gamma_+}
c_1$. Rewrite the function $c_1$ in the form $c_1=g t^{-n}$ and
approximate $g$ by a function $g_{\ve}$ analogously to the
corresponding function in the proof of Proposition \ref{p3}. Recall
that the operator $T(g)$ is invertible and if $\ve$ is small enough,
then
 $$
\boldsymbol\sigma(c_{\ve})=\boldsymbol\sigma(g_{\ve}),
 $$
where $c_{\ve}=\vp_{\gamma_+}g_{\ve} t^n$. Indeed, it follows that
the product of the "+''-factors of the Wiener--Hopf factorization of
the factors in the representation of $c_{\ve}$ is the "+''-factor in
the Wiener--Hopf factorization of $c_{\ve}$ (use \cite[Corollary
5.15]{LiSp:1987}). It remains to show that for sufficiently small
$\ve>0$ one has $\boldsymbol\sigma(c)=\boldsymbol\sigma(c_{\ve})$
but this can be done analogously to the proof of Proposition
\ref{p3}.
 \rbx

\def\cprime{$'$} \def\cprime{$'$}

\end{document}